\documentclass[preprint,12pt]{elsart}

\usepackage{graphicx}
\usepackage{caption}
\usepackage{subcaption}
\usepackage{amssymb}
\usepackage{lineno}

\usepackage{amsmath}
\usepackage{amsfonts}
\usepackage{amssymb}
\usepackage{amscd}
\usepackage{bbm}
\usepackage{bm}
\usepackage{subfig}
\usepackage{epstopdf,epsfig}
\newtheorem{remark}{Remark}

\newcommand{\RR}{\mathbbm{R}}

\newcommand{\dt}{{\Delta t}}

\def\v0{ {{\bf 0}} }

\def\vtau{\mbox{\boldmath $\tau $}}
\journal{Journal of Scientific Computing}

\begin{document}

\begin{frontmatter}



\title{Error Inhibiting Block One-Step Schemes for Ordinary Differential Equations}
%\tnoteref{t1}}
%\tnotetext[t1]{This research was supported by the ISRAEL SCIENCE
%FOUNDATION (grant No. 1364/04) and the UNITED STATES-ISRAEL
%BINATIONAL SCIENCE FOUNDATION (grant No. 2004099).}
%\author[tau]{S. Abarbanel} \ead{saul@post.tau.ac.il}
%\author[tau]{A.�Ditkowski\corref{cor}} \ead{adid@post.tau.ac.il}
%\cortext[cor]{Corresponding author}
%\address[tau]{School of Mathematical Sciences, Tel Aviv University, Tel Aviv 69978, Israel}
%
%\author{}
%
%\address{}

\author{A. Ditkowski}
\address{School of Mathematical Sciences, Tel Aviv University, Tel Aviv 69978, Israel}
\ead{adid@post.tau.ac.il}

\author{S. Gottlieb}
\address{Department of Mathematics, University of Massachusetts, Dartmouth, 285 Old Westport Road, North Dartmouth, MA 02747}
\ead{sgottlieb@umassd.edu}

\begin{abstract}
The commonly used one step methods and linear multi-step methods all have a global error  that is 
of the same order as the local truncation error (as defined in %\cite{AllenIsaacson,IsaacsonKeller,Sewell}
\cite{gustafsson1995time,quarteroni2010numerical,AllenIsaacson,IsaacsonKeller,Sewell}). 
In fact, this is true of the entire class of general linear methods.
In practice, this means that the order of the method is typically defined solely by order conditions which are
derived by studying the local truncation error. In this work we investigate the interplay between the local
 truncation error and the global error, and develop a methodology which defines the construction of explicit {\em error
 inhibiting} block one-step methods (alternatively written as explicit general linear methods \cite{butcher1993a}).
 These  {\em error inhibiting schemes} are constructed so that the accumulation of the local truncation error over time
 is controlled,  which results in a  global error that is one order higher than the local truncation error. 
 In this work, we delineate how to carefully choose the coefficient matrices  so that the  growth of the local truncation error is inhibited.
 We then use this theoretical understanding to construct several methods that have higher order global error 
 than local truncation error, and demonstrate their  enhanced order of accuracy on test cases.
 These methods demonstrate that the error inhibiting concept is realizable. Future work will further develop new 
 error inhibiting methods and will analyze the computational efficiency and linear stability properties of these  methods.
% these {\em error inhibiting schemes} by analyzing the relationship between the local truncation error and the global error. 
% We show that by carefully choosing the coefficient matrices we can ensure that it had only one eigenvalue which is amplified
% over repeated steps, and that the growth of this eigenvalue is mitigated  by the fact that it corresponds to the
%  eigenvector that is a factor of $\Delta t$ times the local truncation error.

%\indent EIS
\end{abstract}

\begin{keyword}
ODE solvers, General linear methods, One-step methods, Global error, local truncation error, Error inhibiting  schemes.
\end{keyword}

\end{frontmatter}

%\linenumbers

\section{Introduction}\label{section:intro}
When solving an ordinary differential equation (ODE) of the form
\begin{eqnarray}\label{ODE}
& & u_t =   F(t,u)  \;,\;\;\;\;\;  t \ge 0  \\
& & u(t=0) =u_0 \; \nonumber
\end{eqnarray}
%(where $F(t,u)$ is continuous with respect to $t$ and Lipschitz-continuous with respect to $u$) 
One can evolve the solution forward in time using the first order forward Euler method
\[ v_{n+1} = v_n + \Delta t F(t_n, v_n) \; .\] 
To obtain a more accurate solution, one can use  methods with multiple steps:
\begin{equation}\label{multistep}
v_{n+1} \,=  \sum_{j=1}^s a_j\,  v_{n+1-j}  + \, \Delta t  \sum_{j=0}^s b_j F(t_{n+1-j}, v_{n+1-j} ),
\end{equation}
known as linear multistep methods \cite{butcher2008numerical}. Alternatively, one can use multiple stages, 
such as Runge--Kutta methods \cite{butcher2008numerical}:
\begin{eqnarray*}\label{RK_Butcher}
y_i  & = & F(v_n+ \sum_{j=1}^m a_{ij} y^{(j)}, t_{n} + \Delta t \sum_{j=1}^m a_{ij}) \; \; \;  \mbox{for} \; \; j=1,...,m\\
v_{n+1}& =& v_n + \Delta t \sum_{j=1}^m b_j y_j .\\
\end{eqnarray*}
The class of general linear methods described in \cite{butcher1993a,JackiewiczBook}
combines the use of multiple steps and stages, constructing methods of the form:
\begin{eqnarray} \label{GLM}
y_i & = &  \sum_{j=1}^s \tilde{U}_{ij} v_{n} +  \Delta t \sum_{j=1}^m \tilde{A}_{ij} f(y_j) \nonumber \\
v^i_{n+1} & =& \sum_{j=1}^s \tilde{V}_{ij} v^i_{n} +  \Delta t \sum_{j=1}^m \tilde{B}_{ij} f(y_j) \; .
\end{eqnarray}
The inclusion of multiple derivatives, such as Taylor series methods \cite{butcher2008numerical}, 
\[ v_{n+1}= v_n + \Delta t F(t_n, v_n) + \frac{\Delta t^2 }{2} F'(t_n,v_n)   + \frac{\Delta t^3 }{3!} F''(v^n), \]
is another possibility, and multiple stages and derivatives have also been developed and used successfully
 \cite{shintani1971one}, \cite{shintani1972explicit},  \cite{kastlunger1972runge}, \cite{kastlunger1972turan}, 
 \cite{chan2010explicit}.

For time-dependent problems  the global error, which is the difference between the 
numerical and exact solution at any given time $t_n= n \Delta t$:
\[ E_n = v_n - u(t_n) ,\]
 depends on  the local truncation error  which, roughly speaking,
 is the accuracy over one time step. % (or every point in space/time).  
 In our case we define the local truncation error  as the error of the method over one time-step, 
 normalized by $\Delta t$.  For example, the local truncation error for Euler's method is (following %\cite{AllenIsaacson,IsaacsonKeller,Sewell}
  \cite{gustafsson1995time,quarteroni2010numerical,AllenIsaacson,IsaacsonKeller,Sewell})
 \[ \tau = \frac{u(t_{n+1})  - u(t_{n})  - \Delta t F(t_n, u(t_{n}) )}{ \Delta t} \approx O(\Delta t). \]
 (To avoid confusion it is important to note that sometimes the truncation error is defined a little differently than we define it above  
 and is not normalized by $\dt$).
 
%(To avoid confusion it is important to note that sometimes the truncation error is defined a little differently than we define it below,  
% and is not normalized by $\dt$, in which case the global error is one order lower than the truncation error).
 
 A well known theoretical result, known as  the
Lax-Richtmeyer  equivalence theorem (see e.g. \cite{lax1956survey}, \cite{gustafsson1995time}, \cite{quarteroni2010numerical}) 
states that if the numerical scheme is stable then the global error is at least of the same order as the local truncation error. 
In all the schemes for numerically solving ordinary differential equations (ODEs) 
that we are familiar with from the literature, the global errors are indeed of the same order as their local truncation errors\footnote{In the case where the truncation error is defined without the  $\dt$ normalization the global error is one order lower than the truncation error.}. 
This is common to other fields in numerical mathematics, such as for finite difference schemes for partial differential equations (PDEs),
see e.g. \cite{gustafsson1995time,quarteroni2010numerical}. It was recently demonstrated, however, that finite difference schemes for PDEs can be constructed such that their convergence rates, or the order of their global errors, are higher than the order of the  truncation errors \cite{ditkowski2015high}. 
In this work we adopt and adapt the ideas presented in \cite{ditkowski2015high} to show that it is possible to construct numerical methods for ODEs
where the the global error is {\em one order higher} than the  local truncation error. As we discuss below, 
these schemes achieve this higher order by inhibiting the lowest order term in the local error from accumulating over time, 
and so we name them {\em  Error Inhibiting Schemes}.

Following an idea in \cite{Rosser1967}, an interesting paper by Shampine and Watt in 1969 \cite{ShampineWatt1969} 
describes a class of implicit one-step methods that obtain a block of
$s$ new step values at each step. These methods take $s$ initial step values and generate the next $s$ step values, and so on, 
all in one step. These methods are in fact explicit block one-step methods, and can be written as general linear methods of the form  \eqref{GLM} above.
Inspired by this form, we construct explicit block one-step methods which are in the form \eqref{GLM}, but where the 
matrix $\tilde{U}$ is an identity matrix, and the matrix $\tilde{A}$ is all zeros; these are known as Type 3 methods in \cite{butcher1993a}.
The major feature of our methods is that in addition to satisfying the appropriate order conditions  listed in \cite{butcher1993a},
they have a special structure  that mitigates the accumulation of  the truncation error,  so we obtain a global error that is one order  
{\em higher}   than predicted by the order conditions  in \cite{butcher1993a}, which describe the local truncation error.

% 
% In Section \ref{section:background} we motivate our approach by describing how
% typical multistep methods can be written and analyzed as block one-step methods: these 
% methods obtain a block of $s$ new step values at each step. We show how this form allows us to
% precisely describe the growth of the local truncation error over the time-evolution. We then exploit this understanding 
% to develop explicit error inhibiting block one-step methods that produce higher order global errors than possible for 
%typical multistep methods in Section \ref{section:EIS}.
% In Section  \ref{section:examples} we present some methods developed according to the theory in 
%Section \ref{section:EIS} and we test these methods on several numerical examples to demonstrate 
%  that the order of convergence is indeed one order higher than the local truncation error.
%  We also show that, in contrast to our error inhibiting Type 3 method, a typical Type 3 method 
%  developed by Butcher in \cite{butcher1993a} does not satisfy the
%  critical condition for a method to be error inhibiting and therefore produces a global error that is of the same
%  order as the local truncation error. Finally, we present our conclusions in Section \ref{conclusions}, and suggest
%that  further investigation  of error inhibiting methods shall include the analysis of their linear stability properties,
%storage implications, and computational efficiency.

 In Section \ref{section:background} we motivate our approach by describing how
 typical multistep methods can be written and analyzed as block one-step methods: these 
 methods obtain a block of $s$ new step values at each step. We show how this form allows us to
 precisely describe the growth of the  error over the time-evolution. 
 In Section \ref{section:EIS} we then exploit this understanding 
 to develop explicit error inhibiting block one-step methods that produce higher order global errors than possible for 
typical multistep methods.
 In Section  \ref{section:examples} we present some methods developed according to the theory in 
Section \ref{section:EIS} and we test these methods on several numerical examples to demonstrate 
  that the order of convergence is indeed one order higher than the local truncation error.
  We also show that, in contrast to our error inhibiting Type 3 method, a typical Type 3 method 
  developed by Butcher in \cite{butcher1993a} does not satisfy the
  critical condition for a method to be error inhibiting and therefore produces a global error that is of the same
  order as the local truncation error. Finally, we present our conclusions in Section \ref{conclusions}, and suggest
that  further investigation  of error inhibiting methods shall include the analysis of their linear stability properties,
storage implications, and computational efficiency.

\section{Motivation}\label{section:background}
In this section we present the analysis of explicit multistep methods in a block one-step form for a simple linear problem.
In this familiar setting we define the local truncation error, the global error, and the solution operator that connects them.
We also discuss the stability of a  method of this form.
We limit our analysis to the linear case so that we can clearly observe the process by which the 
solution operator interacts with the local truncation error, and results in a global error that is of the same order as the local
truncation error. Although we are dealing for the moment 
with standard multistep methods, this will set the stage for the construction and analysis of error inhibiting block one-step
methods.

In order to illustrate the main idea we start with a linear ordinary differential equation (ODE) 
\begin{eqnarray}\label{1.10}
& & u_t =   f(t) \; u  \;,\;\;\;\;\;  t \ge 0  \\
& & u(t=0) =u_0 \; \nonumber
\end{eqnarray}
where $f(t)<M \,  , \;  \forall t \ge 0$ and $f(t)$ is analytic. 

An  $s$-step explicit  multistep method applied to \eqref{1.10}  takes the form
\begin{equation}\label{Standard_multistep_1}
v_{n+s} = \sum_{j=0}^{s-1} a_j\, v_{n+j} \, + \, \dt \sum_{j=0}^{s-1} b_j F(t_{n+j},  \, v_{n+j} )
= \sum_{j=0}^{s-1} a_j\, v_{n+j} \, + \, \dt \sum_{j=0}^{s-1} b_j f(t_{n+j})  \, v_{n+j} 
\end{equation}
where the time domain is discretized by the sequence $t_n = n \, \dt$, and 
$v_n$ denotes the numerical approximation of $u(t_n)$. The method \eqref{Standard_multistep_1} is defined by its coefficients
$ \{ a_j \}_{j=0}^{s-1}$ and $\{b_j \}_{j=0}^{s-1}$, which are constant values. 
%For this scheme to be explicit, $b_s=0$ (and $a_s \neq 1$). For convenience, we define $a_s=1$.

Following \cite{gustafsson1995time} we rewrite the method \eqref{Standard_multistep_1} in its block form. To do this,  
we first introduce the exact solution vector
\begin{equation}\label{multistep_u}
U_n = \left( u(t_{n+s-1}), \ldots, u(t_n) \right )^T
\end{equation}
and similarly, the numerical solution vector is
\begin{equation}\label{multistep_v}
V_n = \left( v_{n+s-1}, \ldots, v_n\right )^T .
\end{equation}
Now  \eqref{Standard_multistep_1} can be written in block form so that it looks like a one step scheme
\begin{equation} \label{Standard_multistep_2}
V_{n+1} = Q_n V_n 
\end{equation}
where
\begin{equation}\label{Q_def}
Q_n = \left(  \begin{array}{cccc}
a_{s-1}+\dt  b_{s-1} f(t_{n+s-1}) \; \; & \; \; a_{s-2}+\dt b_{s-2}f(t_{n+s-2}) \;  \; & \dots & \; \; a_{0}+\dt b_{0} f(t_{n})\\
I \\
& \ddots \\
& & I &0
\end{array}  \right )  .
\end{equation}

From repeated applications of equation \eqref{Standard_multistep_2} we observe that the numerical solution vector  $V_n$ at any time $t_n$ 
can be related to $V_{\nu}$ for any  previous time  $t_\nu$ 
\begin{equation}\label{S_h_def}
V_n = S_\dt \left( t_n, t_{\nu}\right ) V_{\nu}\;,\;\;\;{\nu} \le n
\end{equation}
where $S_\dt$ 
 is the discrete solution operator. This operator can be expressed explicitly by
\begin{equation}\label{S_h_def_2}
S_\dt \left( t_n, t_{\nu}\right ) =  Q_{n-1} \ldots Q_{\nu+1} Q_{\nu}\;,\;\;\;S_\dt \left( t_n, t_n\right )  = I.
%S_\dt \left( t_n, t_{\nu}\right ) = \prod_{\mu=1}^{n-\nu} Q_{n-\mu} \;,\;\;\;S_\dt \left( t_n, t_n\right )  = I.
\end{equation}
For simplicity  we can define this by
\begin{equation}\label{S_h_def_2.5}
\prod_{\mu=\nu}^{n-1} Q_{\mu}\equiv  Q_{n-1} \ldots Q_{\nu+1} Q_{\nu}\;,\;\;\;\prod_{\mu=n}^{n-1}Q_{\mu} \equiv I.
\end{equation}
Note that if each matrix $Q_\mu$ is independent of $\mu$ (in other words, in the constant coefficient case where $f$ is independent of $t$),
we simply have a product of matrices $Q$, and the discrete solution operator becomes 
\begin{equation}\label{S_h_def_3}
S_\dt \left( t_n, t_{\nu}\right ) = Q^{n-\nu}.
\end{equation}

The behavior of a method depends in large part on the accuracy of its solution operator. We begin by defining the local truncation error
 as the error of the method over one time-step, normalized by $\dt$: \newline
\noindent{\bf Definition 1}: \label{def:truncation_error} 
{\em The local  truncation error $\vtau_n$ is given by %\cite{Sewell,AllenIsaacson,IsaacsonKeller}
\cite{gustafsson1995time,quarteroni2010numerical,AllenIsaacson,IsaacsonKeller,Sewell}
\begin{equation} \label{truncation_error_1}
\dt \, \vtau_n \, = \, U_{n+1} - Q_n U_n
\end{equation}
}

Note that in the case of the standard multistep method, where $Q_n$ is given by the matrix
 \eqref{Q_def}, the truncation error has only one non-zero element:
\begin{equation} \label{truncation_error_2}
\vtau_n \, = \, \left (  \tau_n, 0, \ldots, 0\right )^T .
\end{equation}

The error that we are most interested in  is the difference between the exact error vector and the numerical error vector
 at time $t_n$, 
\begin{equation}  \label{error}
E_n = U_n-V_n \; ,
\end{equation}
known as the global error.
At the initial time, we have the error $E_0$ which is based on the starting values a method of this sort requires:
the values $v_j$, $j=0, \ldots, s-1$  that are prescribed or somehow computed. 
Typically, $v_0$ is the initial condition defined  in \eqref{ODE} and $v_j$, $j=1, \ldots, s-1$ are computed to sufficient accuracy
using some other numerical scheme. Thus, the value of $E_0$ is assumed to be as small as needed.

The evolution of the global error \eqref{error} depends on  the local truncation error defined by \eqref{truncation_error_1}
and the discrete solution operator given in \eqref{Standard_multistep_2}:
\begin{equation} \label{error_equation_1}
E_{n+1} \, = \,  Q_n E_n     \,+\,      \dt \, \vtau_n \;.
\end{equation}
Unraveling this equality all the way back to $E_0$ gives
\begin{equation} \label{error_equation_2}
E_{n} \, = \,  S_\dt \left( t_n, 0\right ) E_0  \,+\,  \dt \, \sum_{\nu=0}^{n-1}  S_\dt \left( t_n, t_{\nu+1}\right ) \vtau_{\nu} \;,
\end{equation}
%\begin{equation} \label{error_equation_2a}
%E_{n} \, = \,   \prod_{\mu=1}^{n} Q_{n-\mu}  E_0  \,+\,  \dt \, \sum_{\nu=0}^{n-1}  
%\left( \prod_{\mu=1}^{n-\nu -1} Q_{n-\mu-1}  \right)  \vtau_{\nu} \;.
%\end{equation}
or, equivalently
%\begin{equation} \label{error_equation_2}
%E_{n} \, = \,  S_\dt \left( t_n, 0\right ) E_0  \,+\,  \dt \, \sum_{\nu=0}^{n-1}  S_\dt \left( t_n, t_{\nu+1}\right ) \vtau_{\nu} \;.
%\end{equation}
\begin{equation} \label{error_equation_2a}
E_{n} \, = \,   \prod_{\mu=0}^{n-1} Q_{\mu}  E_0  \,+\,  \dt \, \sum_{\nu=0}^{n-1}  
\left( \prod_{\mu=\nu+1}^{n-1} Q_{\mu}  \right)  \vtau_{\nu} \;.
\end{equation}
%\prod_{\mu=\nu}^{n-1} Q_{\mu}

(This  formula is obtained from the discrete version of Duhamel's principle, see Lemma 5.1.1 in  \cite{gustafsson1995time}).

%It is clear from \eqref{error_equation_2} that the behavior of the discrete solution operator $S_\dt ( t_n, t_{\nu+1}) $ must be 
%controlled for this error to converge. This property defines the stability of the method:  \newline
It is clear from \eqref{error_equation_2} that the behavior of the discrete solution operator $S_\dt ( t_n, t_{\nu+1}) $ must be 
controlled for this error to converge. This property defines the stability of the method. Also here we use the stability definition presented in  \cite{gustafsson1995time}, namely: \newline
%The behavior of a scheme, and in particular its stability, is then determined by the behavior of the solution operator: 
 \noindent{\bf Definition 2}:    \label{def:scheme_stability}
{\em The scheme \eqref{Standard_multistep_2} is called stable  if there are constants  $\alpha_s$ and $K_s$, independent of $\dt$,  
 such that for all $0 <\dt  \le \dt_0$
\begin{equation}\label{S_h_stability_1}
\left \| S_\dt \left( t_n, t_{\nu}\right )  \right \| \le K_s e^{\alpha_s\left( t_n - t_{\nu}\right )  }
\end{equation}}

If the scheme is stable, we can use \eqref{S_h_stability_1} and   \eqref{error_equation_2}   to bound the growth of the error:
\begin{equation} \label{error_equation_3}
\left \|  E_{n}   \right\| \, \le \,  K_s \left [ e^{\alpha_s  t_n   } \left \|  E_{0}   \right\|\,+\,  \max_{0 \le \nu \le n-1} \left \|  \vtau_{\nu}   \right\| \phi_h^*(\alpha_s, \, t_n) \right ] \;.
\end{equation}
where
\begin{equation} \label{error_equation_4}
\phi_{\dt}^*(\alpha_s, \, t_n) \,=\, {\dt} \,\sum_{\nu=0}^{n-1} e^{\alpha_s\left( t_n - t_{\nu+1}\right)} \, \approx\, 
\int_0^{t_n} e^{\alpha_s\left( t_n - \zeta\right)}d\zeta
  \,=\, \left \{\begin{array}{lll}
\frac{e^{\alpha_s\, t_n } -1}{\alpha_s} & \hspace{1cm} & \alpha_s \neq 0  \\
t_n && \alpha_s = 0
\end{array} \right. \;.
\end{equation}
Equation \eqref{error_equation_3} means that stability implies convergence:\footnote{
For partial differential equations this result is known as one part of the celebrated 
Lax-Richtmeyer  equivalence theorem. See e.g. \cite{lax1956survey}, \cite{gustafsson1995time}, \cite{quarteroni2010numerical}.}
 if the scheme is stable than the global error is controlled by the local truncation error for any given final time.  
 In the formula above it is clear that the global error must have order at least  as high as the local truncation error, 
 but the possibility of having a higher order global error is left open.

%By taking the norm of  \eqref{error_equation_2} and using \eqref{S_h_stability_1} and the  inequality, $\left \|  E_{n}   \right\|$ can be bounded by:
The first Dahlquist barrier \cite{hairer2000solving,butcher2008numerical} states that any explicit $s$ step linear
multistep method can be of order $p$ no higher than $s$. It is the common experience that methods
have global error of the same order as the local truncation error. These two together greatly limit the accuracy of the
methods we can derive. 

% rem \ref{remark:Adams-Bashforth_rem} 
%{\bf Remark 1}: \label{remark:Adams-Bashforth_rem} {\em In an Adams-Bashforth scheme the entry in the first row and first column in the term  
%$S_\dt \left( t_n, t_{\nu}\right ) = \prod_{\mu=\nu}^{n-1} Q_{\mu}$  is equal to $1+O(\dt)$. Therefore the error, 
%due to the accumulation of the contributions from the truncation errors,  becomes:
%\begin{equation} \label{AdamsBashfortherro_1}
%{\bf e}_{n+s} \, = \,   \dt \, \sum_{\nu=0}^{n-1}  \left(   1+O(\dt)  \right ) \tau_{\nu} \;
%\end{equation}
%which is approximately the  average value  of $\tau_{\nu}$  over ${\nu=0,..,n-1} $.
%This  suggests that we may need to look outside the family of linear multistep methods to attain a higher order
%global error. }

\begin{remark}\label{remark:Adams-Bashforth_rem} 
%{\em 
In an Adams-Bashforth scheme the entry in the first row and first column in the term  
$S_\dt \left( t_n, t_{\nu}\right ) = \prod_{\mu=\nu}^{n-1} Q_{\mu}$  is equal to $1+O(\dt)$. Therefore the error, 
due to the accumulation of the contributions from the truncation errors,  becomes:
\begin{equation} \label{AdamsBashfortherro_1}
{\bf e}_{n+s} \, = \,   \dt \, \sum_{\nu=0}^{n-1}  \left(   1+O(\dt)  \right ) \tau_{\nu} \;
\end{equation}
which is approximately the  average value  of $\tau_{\nu}$  over ${\nu=0,..,n-1} $.
This  suggests that we may need to look outside the family of linear multistep methods to attain a higher order
global error.
% }
\end{remark}

The analysis in this section suggests that if the operator $Q_n$ is properly constructed,
the growth of the global  error described in Equation \eqref{error_equation_2a} may be controlled
 through the properties of the operator $Q_n$ 
and its relationship with the  local truncation error $\vtau_{n}$. However, as implied by the example of the 
Adams-Bashforth scheme above, we need to construct methods where the operator $Q_n$ is not limited
by the structure in this section. In the next section we present the construction of block one-step
methods that are error inhibiting.
The class of methods described  by this block one-step  structure is very broad: while all classical multistep 
 methods can be written in this block form, not  every such block one-step method can be 
 written as a classical multistep method. Thus, we rely on the discussion in this section with one main change:
 the structure of the  operator $Q_n$.

%%%%%%%%%%%%%%%%%%%%%%%%%%%%%%%%%%%%%%%%%%%

\section{An Error Inhibiting Methodology} \label{section:EIS}
In Section \ref{section:background} we rewrote  explicit  linear   multistep methods in a block one-step form, 
and expressed the relationship between its local and global error. We observed that the growth of the local errors 
is driven by the behavior of the discrete solution operator $Q_n$, and in particular its interaction with the 
local truncation error. Using this insight
%This analysis raises the question of whether one can construct  variants of the multistep method 
%that are not limited by an order barrier of the type shown by Dahlquist for multistep methods. 
%{\bf Sigal will think about putting it in general???}
we show in this section that it is possible to construct such explicit block one-step methods 
(which are also known as Type 3 DIMSIM methods in \cite{butcher1993a})
that {\em inhibit} the growth of the truncation error so that the global error \eqref{error} gains an order of 
accuracy over the local truncation error  \eqref{truncation_error_1}. 

We begin in Section \ref{section:EIS_CC_1} by describing the construction and analysis of 
error inhibiting block one-step schemes for the case of  linear constant coefficient equations.
We then  show that this approach yields methods that are also error inhibiting for variable coefficient 
linear equations in Section \ref{section:EIS_VC_1}
and nonlinear equations in Section \ref{section:EIS_NL_1}.

\subsection{Error inhibiting schemes for linear constant coefficient equations} \label{section:EIS_CC_1}

Given a linear ordinary differential equation with  constant coefficients:
\begin{eqnarray}\label{EIS_10}
& & u_t=  f \cdot u   \;,\;\;\; \mbox{for} \; \; \; \; t \ge 0,  \\
& & u(t=0) = u_0 \nonumber 
\end{eqnarray}
where $f \in \RR$. 
We define a vector of length $s$ that contains the exact solution of \eqref{EIS_10} 
at times  $\left( t_n+ j \Delta t/s \right)$ for $j=0,\ldots,s-1$
\begin{equation}\label{multistep_u}
U_n = \left( u(t_{n+(s-1)/s}), \ldots,  u(t_{n+1/s}),  u(t_n) \right )^T ,
\end{equation}
and the corresponding vector of numerical approximations
\begin{equation}\label{multistep_v}
V_n = \left( v_{n+(s-1)/s}, \ldots,  v_{n+1/s},  v_n\right )^T.
\end{equation}

Note that although we are assuming that the solution $u$ at any given time is a scalar, 
this entire discussion easily generalizes to the case where $u$ is a vector, with only 
some cumbersome notation needed. Thus without loss of generality we continue the discussion 
with scalar notation.

%{\bf Remark}: {\em The notation above emphasizes that this scheme uses $s$ terms for generating the next $s$ terms, 
%unlike the  explicit linear multistep methods in the section above which use $s$ terms to generate one term. 
%To match with the notation in Section 2 above,  we can replace $\Delta t' = s \Delta t $ thus  defining this scheme 
%on integer grid points. }

\begin{remark}\label{remark:s_explanation_rem} 
The notation above emphasizes that this scheme uses $s$ terms for generating the next $s$ terms, 
unlike the  explicit linear multistep methods in the section above which use $s$ terms to generate one term. 
To match with the notation in Section 2 above,  we can replace $\Delta t' = s \Delta t $ thus  defining this scheme 
on integer grid points.
\end{remark}

We define the  block one-step method for the linear constant coefficient problem \eqref{EIS_10}  
\begin{equation} \label{EIS_multistep_2}
V_{n+1} = Q V_n
\end{equation}
where
\begin{equation}\label{EIS_Q_def}
Q = A + \Delta t B f 
\end{equation}
and  $A, B \in \RR^{s \times s}$. Unlike in the case of classical multistep methods, here we do not 
restrict the structure of the matrices $A$ and $B$. Thus, any multistep method of the form \eqref{Standard_multistep_1} can be written in this form
(as we saw above), but not every method of the form 
\eqref{EIS_Q_def} can be written as a multistep method.
In fact, this methods is a general linear method of the DIMSIM form \eqref{GLM} 
with $\tilde{A}$ is all zeroes, $\tilde{U} $ is the identity matrix,
$\tilde{V}=A$, and  $\tilde{B}= B $.
This particular  formulation is, as we mentioned above, 
called a Type 3 DIMSIM in Butcher's 1993 paper \cite{butcher1993a}.

At any time $t_n$, we know that $ \; u(t_{n}+\dt) \,= \, u(t_n)+O(\dt)$, so that for the numerical solution $V_n$ to 
converge to the analytic solution $U_n$  one of the eigenvalues of $Q$ must be equal to $1+O(\dt)$,
 and its eigenvector must have the form:
\begin{equation} \label{EIS_20}
\left( 1+O(\dt), \ldots, 1+O(\dt)\right)^T \; .
\end{equation}
The structure of the eigensystem of $A$, which is the leading part of $Q$, is critical to the stability of the scheme
and the dynamics of the error.

Suppose $A$ is constructed such that:
\begin{enumerate}  %% to rewrite later Sum of rows=1 ????
\item ${\bf rank} (A)=1$.
\item  Its non-zero eigenvalue is equal to one and its corresponding eigenvector is 
$  \left( 1, \ldots, 1\right)^T $
\item $A$ can be diagonalized. %Furthermore we assume that there is a constant $K_A$ ????
\end{enumerate}
Property (2) assures that the method produces the exact solution for the case $f=0$.
Now, since the term $ \dt B f $ is only an $O(\dt)$ perturbation to $A$, the matrix $Q$ will have one eigenvalue, 
$z_1=1+O(\dt)$  whose eigenvector has the form
\begin{equation} \label{EIS_40}
\psi_1 = \left( 1+O(\dt), \ldots, 1+O(\dt)\right)^T\;
\end{equation}
and the rest of the eigenvalues satisfy $z_j=O(\dt)$ for $j=2,\ldots,s$. 

Since the  $\|Q\| = 1+O(\dt)$, we can conclude that there exist constants $K_s$ and $\alpha_s$
such that 
\begin{equation}\label{EIS_stability_1}
\left \| S_\dt \left( t_n, t_{\nu}\right )  \right \|  = \left\| Q^{n-\nu} \right \|  \le K_s e^{\alpha_s\left( t_n - t_{\nu}\right )  }
\end{equation}
where $\alpha_s = \|B\| \, |f|$. Therefore, according to Definition 2, the scheme \eqref{EIS_multistep_2} is stable.
By the same argument used above, we can show that the global error will have order that is no less than the order of the
local truncation error.

We now turn to the task of investigating the truncation error, $\vtau_n$. The definition of the
local  truncation error in this  case is still  
\[ \dt \, \vtau_n \, = \, U_{n+1} - Q_n U_n \]
as defined in the previous section in Equation  \eqref{truncation_error_1}.

\begin{remark}\label{remark:truncation_error_rem} 
Since $Q = A + \Delta t B f$ and $u_t=fu$ the local truncation error can be written as
$$
\dt \, \vtau_n \, = \, U_{n+1} - \left( A U_n +\dt B \frac{d\, U_n}{dt}\right) \;.
$$
Therefore $\vtau_n$ does not explicitly depend on $f$. This observation is valid for the variable coefficients and the nonlinear case as well.
\end{remark}

 %The equation for the error is 
 The definition of the error is 
\[E_n = U_n-V_n \; ,\]
 as in Equation \eqref{error}. The evolution of the error is still described by
Equation \eqref{error_equation_2a}
%\[ E_{n} \, = \,   \prod_{\mu=1}^{n} Q_{n-\mu}  E_0  \,+\,  \dt \, \sum_{\nu=0}^{n-1}  
%\left( \prod_{\mu=1}^{n-\nu -1} Q_{n-\mu-1}  \right)  \vtau_{\nu} \;,\]
$$
E_{n} \, = \,   \prod_{\mu=0}^{n-1} Q_{\mu}  E_0  \,+\,  \dt \, \sum_{\nu=0}^{n-1}  
\left( \prod_{\mu=\nu+1}^{n-1} Q_{\mu}  \right)  \vtau_{\nu} \;,
$$
which in the linear constant coefficient case becomes
\begin{equation} \label{error_equation_cc}
%E_{n}  =   \prod_{\mu=1}^{n} Q^{n-\mu}  E_0  \,+\,  \dt \, \sum_{\nu=0}^{n-1}  Q^{n-\mu-1}   \vtau_{\nu} \;.
E_{n}  =   Q^{n}  E_0  \,+\,  \dt \, \sum_{\nu=0}^{n-1}  Q^{n-\nu-1}   \vtau_{\nu} \;.
\end{equation}
The main difference between this case and the linear multistep method in Section \ref{section:background} is that
the structure of $Q$ is different, and that  unlike \eqref{truncation_error_2},  in this case all the entries in $\vtau_n$ are typically non-zero.

Equation \eqref{error_equation_cc} indicates  that there are several sources for the error at the time $t_n$:
\begin{enumerate}
\item {\em The initial error $E_0$  which is the error in the initial condition $V_0$:}  This error is  caused primarily 
by the numerical scheme used to compute the first $s-1$ elements  in $V_0$.  We assume these errors can be made arbitrary small.
The initial value, which is the final element of $V_0$, is taken from the analytic initial condition and is considered to be accurate to machine precision. 
\item {\em The term  $\dt \,\vtau_{n-1},$ which is  the last term in the sum in the right hand side of \eqref{error_equation_cc}}:
This term  is clearly, by definition, of the size $O(\dt) \|\vtau_{n-1}\|.$
\item {\em The summation 
\begin{eqnarray} \label{error_sum}
\dt \, \sum_{\nu=0}^{n-2} Q^{n - \nu -1}  \vtau_{\nu},
\end{eqnarray} which are all the 
rest of the  terms in the sum in the right hand side of \eqref{error_equation_2}}: This is the term we need to bound to control
the growth of the truncation error. 
\end{enumerate}

The terms in the sum \eqref{error_sum} are all comprised of the discrete solution operator $Q$ multiplying the 
local truncation error. This leads us to the major observation that is the key to constructing error inhibiting methods:
{\bf if the local truncation error lives in the subspace of eigenvectors that correspond to the eigenvalues of $O(\dt)$,
then the growth of the truncation error will be inhibited, and the global error will be one order higher than the local
truncation error.}

Recall that  $Q$ has one dominant eigenvalue  that has the form $ 1+O(\dt)$  and all the others are $O(\dt)$. 
Correspondingly, two subspaces can be defined
\[  \Psi_1  = {\rm span} \left \{\psi_1\right \} \; \; \; \; \mbox{and} \; \; \; \; 
\Psi_1^c  = {\rm span} \left \{\psi_2, ..., \psi_s  \right  \}
\]
where $\psi_j$ is the eigenvector associated with each eigenvalue $z_j$. As $\psi_j$ can be normalized, we assume that $\| \psi_j \|=O(1)$. 
It should be noted that while $\Psi_1$ and $\Psi_1^c$  are linearly independent, they are not orthogonal subspaces.
Furthermore, since the   matrix $A$ is  diagonalizable by construction,  its eigenvectors span $\RR^s$. 
Since $ \vtau_{\nu} \in \RR^s$,  it can be written as
\begin{equation}\label{vtau_expansion_1}
\vtau_{\nu} \, = \,  \gamma_1 \psi_1 + \sum_{j=2}^s \gamma_j \psi_j   \;
\end{equation}
where  $\gamma_1 \psi_1 \in \Psi_1 $ and $\sum_{j=2}^s \gamma_j \psi_j  \in \Psi_1^c$. 

Of course, the truncation error $\vtau_{\nu} $  is determined by the entries of $Q$. 
To ensure that the local truncation error is mostly in the space $\Psi_1^c $ of eigenvectors which 
correspond to the eigenvalues of size $O(\dt)$, we  choose the entries of $Q$ (i.e. the entries of $A$ and $B$) 
such that  $\gamma_1 = O(\dt)$, which will mean that 
\begin{equation}\label{construction_criteria_1}
\left \| \gamma_1 \psi_1 \right \|  \, = \,  O(\dt) \left \| \vtau_{\nu} \right \|   \;.
\end{equation}
Using this, we can bound product of the discrete solution operator and the truncation error,
\begin{eqnarray}\label{construction_criteria_1.5}
\left \|  Q \vtau_{\nu} \right \|  & = & \left \| \gamma_1 Q \psi_1 + \sum_{j=2}^s \gamma_j Q \psi_j \right \|  
\leq  \left \| \gamma_1 Q \psi_1 \right \| +\left \|   \sum_{j=2}^s  \gamma_j Q \psi_j \right \| \nonumber \\
& = & \left \| \gamma_1 z_1 \psi_1 \right \| +\left \| \sum_{j=2}^s   \gamma_j z_j \psi_j \right \| 
\leq |z_1| \left \| \gamma_1 \psi_1 \right \| + \max_{j=2, .. . s}  |z_j|  \left \| \sum_{j=2}^s  \gamma_j \psi_j \right \| \nonumber\\
&\leq & |z_1| \left \| \gamma_1 \psi_1 \right \| + \max_{j=2, .. . s}  |z_j|  \left \| \vtau_\nu - \gamma_1 \psi_1  \right \| \nonumber\\
& \leq & \left(1+ O(\dt) \right) O(\dt)  \left \| \vtau_{\nu} \right \|  + O(\dt)  \left \| \vtau_{\nu} \right \|  
= O(\dt)  \left \| \vtau_{\nu} \right \| \nonumber
\end{eqnarray}
where $z_j$ are the eigenvalues of $Q$. 
%Since $|  \gamma_1 | = O(\dt)$  and $\left \|  \psi_1\right \| = O(1)$,
%$\left \| \sum_{j=2}^s  \gamma_j \psi_j \right \| = O(\vtau_{\nu} )$,  see \eqref{vtau_expansion_1}  above. 
Therefore we have 
\begin{equation} \label{construction_criteria_2}
\left \|  Q \vtau_{\nu} \right \|  \leq O(\dt)  \left \| \vtau_{\nu} \right \| \ .
\end{equation}

Whenever the condition \eqref{construction_criteria_2} is satisfied, we can show that the sum \eqref{error_sum} above is bounded:
\begin{eqnarray} \label{construction_criteria_2.5}
\left\| \dt \, \sum_{\nu=0}^{n-2} Q^{n - {\nu-1}}  \vtau_{\nu}  \right\| &=& 
\dt \left\|  \sum_{\nu=0}^{n-2} Q^{n - {\nu-1}}  \vtau_{\nu}  \right\| \nonumber  \leq 
\dt  \sum_{\nu=0}^{n-2} \left\| Q^{n-\nu-2} \right\| \left\| Q \vtau_{\nu}  \right\| \nonumber  \\
& \leq & \dt \sum_{\nu=0}^{n-2} \left\| Q \right\|^{n-\nu-2} O(\dt)   \|\vtau_\nu \| \nonumber  \\
& \leq & \dt   \left( \max_{\nu=0,...,n-2} \left \| \vtau_{\nu} \right \|  \right)  \sum_{\nu=0}^{n-2} (1+c \dt)^{n-\nu-2} O(\dt)  \nonumber   \\
& \leq & \dt   \left( \max_{\nu=0,...,n-2} \left \| \vtau_{\nu} \right \|  \right)  \sum_{\nu=0}^{n-2}  \left[ e^{c \dt} \left( 1+O(\dt^2) \right) \right ]^{n-\nu-2} O(\dt)  \nonumber   \\
& \leq & \dt   \left( \max_{\nu=0,...,n-2} \left \| \vtau_{\nu} \right \|  \right)  \sum_{\nu=0}^{n-2}  \left[ e^{c (t_{n-2}-t_{\nu})} \left( 1+O(\dt) \right) \right ] O(\dt)  \nonumber   \\
& \leq & O(\dt)   \left( \max_{\nu=0,...,n-2} \left \| \vtau_{\nu} \right \|  \right) \phi_{\dt}^*(c, \, T) . 
% O(\dt) \max_{\nu=0,...,n-2} \left \| \vtau_{\nu} \right \|   \;.
\end{eqnarray}
(Recall \eqref{error_equation_4} for the definition of  of $\phi_{\dt}^*(c, \, T)$.)

%\begin{eqnarray*} 
%\left\| \dt \, \sum_{\nu=0}^{n-2} Q^{n - {\nu-1}}  \vtau_{\nu}  \right\| &=& 
%\dt \left\|  \sum_{\nu=0}^{n-2} Q^{n - {\nu-1}}  \vtau_{\nu}  \right\| \leq 
%\dt  \sum_{\nu=0}^{n-2} \left\| Q^{n-\nu-2} \right\| \left\| Q \vtau_{\nu}  \right\| \\
%& \leq & \dt \sum_{\nu=0}^{n-2} \left\| Q \right\|^{n-\nu-2} O(\dt)   \|\tau_\nu \|  \\
%& \leq & \dt   \left( \max_{\nu=0,...,n-2} \left \| \vtau_{\nu} \right \|  \right)  \sum_{\nu=0}^{n-2} (1+c \dt)^{n-\nu-2} O(\dt)    \\
%& \leq & \dt   \left( \max_{\nu=0,...,n-2} \left \| \vtau_{\nu} \right \|  \right)  \left( \frac{(1+c \dt)^{n-1} -1}{1+c\dt - 1} \right) O(\dt)    \\
%& \leq & O(\dt)  \max_{\nu=0,...,n-2}  \left \| \vtau_{\nu} \right \|    e^{c n \dt}  \leq  O(\dt)  \max_{\nu=0,...,n-2}  \left \| \vtau_{\nu} \right \|    e^{c T}.  \\
%% O(\dt) \max_{\nu=0,...,n-2} \left \| \vtau_{\nu} \right \|   \;.
%\end{eqnarray*}
In the final equation, $T$ is the final time, and the term $\phi_{\dt}^*(c, \, T)$
%$e^{c T}$ 
is therefore a constant. Thus we have the bound
\begin{equation} \label{construction_criteria_3}
\left\| \dt \, \sum_{\nu=0}^{n-2} Q^{n - {\nu-1}}  \vtau_{\nu}  \right\| \leq 
O(\dt)  \max_{\nu=0,...,n-2}  \left \| \vtau_{\nu} \right \|  .
\end{equation}
Putting this all together  into \eqref{error_equation_cc}, we obtain
\begin{equation}\label{EIS_error_estimate_1}
\left \|  E_n \right \|  \, = \, O(\dt) \max_{\nu=0,...,n-1} \left \| \vtau_{\nu} \right \|   \;.
\end{equation}
Thus, if the coefficients of $A$ and $B$ are chosen so that we can 
control the size of $\left \|  Q \vtau_{\nu} \right \|$ in \eqref{construction_criteria_2}, we 
can obtain a scheme that inhibits the growth of the local truncation error, 
so that the global error is  one order more accurate than its truncation error. 

%Since, by \eqref{construction_criteria_1} -- \eqref{construction_criteria_3}, 
%such schemes inhibit the error from  accumulating in time, we name them {\em error inhibiting schemes (EIS)}.

\subsection{Linear variable-coefficient equations} \label{section:EIS_VC_1}
In the previous section we showed how to construct an error inhibiting method by choosing the coefficients in  $A$ and $B$ 
so that the local truncation error lives (mostly) in the space that is spanned by the eigenvectors corresponding to eigenvalues 
that are of $O(\dt)$.  In this section we show that under the same criteria as above, these methods are also error inhibiting when applied to
a {\em variable  coefficient} linear ordinary differential equation:
\begin{eqnarray}\label{EIS_vc_10}
&& u_t=  f(t)  u   \;\;,\;\;\;\;\;  t \ge 0 \nonumber \\
&& u(t=0) = u_0 \;
\end{eqnarray}
where $f(t)$ assumed to be analytic or as smooth as needed, and bounded. 
In this case the scheme is given by a time-dependent evolution operator $Q_n$ which
may change each time-step:
\begin{equation} \label{EIS_multistep_vc_2}
V_{n+1} = Q_n V_n
\end{equation}
where
\begin{equation}\label{EIS_Q_vc_def_1}
Q _n= A + \dt B \,
\left (  \begin{array}{cccccc}
f \left( t_{n+(s-1)/s} \right ) \\
& f \left( t_{n+(s-2)/s} \right ) \\
&& \ddots \\
&&& \;\;\;\; f \left( t_{n} \right )
\end{array} \right )
\end{equation}
and the matrices $A$ and $B$ are the same as described above for the constant coefficient scheme.

Since $f(t)$ is an analytic function, $Q_n$ can be written as
\begin{equation}\label{EIS_Q_vc_def_2}
Q _n= A + \dt B f(t_n) + \dt^2 B f'(t_n)
\left (  \begin{array}{cccccc}
 \left( {(s-1)/s} \right ) \\
&  \left( {(s-2)/s} \right ) \\
&& \ddots \\
&&& \;\;\;\;\; 0
\end{array} \right ) + O(\dt^3)
\end{equation}
We can also say then that
\begin{equation}\label{EIS_Q_vc_def_3}
Q_{n} = A + \dt B f(t_{n}) + O(\dt^2) B f'(t_{n})  =
%A + \dt B f(t_{n}) +  O( \dt^2)  =
 \tilde{Q}_{n} + O( \dt^2) .
\end{equation}
%\[Q_{n-\mu} = A + \dt B f(t_{n-\mu}) + \dt^2 B f'(t_{n-\mu}) + O(\dt^2) =
%A + \dt B f(t_{n}) + \mu O( \dt^2)  \]
Each $\tilde{Q}_{n}$ has the same structure as $Q$ in the constant coefficient case. In particular 
\begin{equation}\label{EIS_Q_vc_def_4}
\| \tilde{Q}_{n} \| = \left(1+O(\dt)\right )\le 1+c \dt,\; \;\;\; \forall n\;.
\end{equation}
Furthermore, as was pointed out in Remark \ref{remark:truncation_error_rem}, since the  local truncation error $\vtau_n$ 
does not depend explicitly on $f(t)$ at any time $t_n$, we can write  $\vtau_n$   as a linear combination of the eigenvectors of 
$A$ that correspond to the zero eigenvalues. Thus,  $\vtau_n$ lives (mostly) in the space that is spanned by the eigenvectors 
of any matrix $\tilde{Q}_{n}$ corresponding to eigenvalues  that are of $O(\dt)$. 
We can then follow the  same analysis as in \eqref{construction_criteria_1}--\eqref{construction_criteria_1.5},
to obtain the bound
\begin{equation}\label{EIS_Q_vc_def_5}
\| \tilde{Q}_{n+1} \vtau_n\| = O(\dt)\| \vtau_n\|  ,\; \;\;\; \forall n \;.
\end{equation}

In this case, Equation \eqref{error_equation_2}  takes the modified form (for $n \geq 1$)
\begin{eqnarray} \label{error_equation_VC_1}
E_{n} & = &  \prod_{\mu=0}^{n-1} Q_{\mu}  E_0  \,+\,  \dt \, \sum_{\nu=0}^{n-1}  
\left( \prod_{\mu=\nu+1}^{n-1} Q_{\mu}  \right)  \vtau_{\nu} \nonumber \\
& = &  \prod_{\mu=0}^{n-1} Q_{\mu} E_0  \,+\,  
\dt \, \sum_{\nu=0}^{n-2}  \prod_{\mu=\nu+1}^{n-1}
 \left ( \tilde{Q}_{\mu} +  O(\dt^2) \right)  \vtau_{\nu}  +
 \dt \vtau_{n-1}
 \nonumber  
 \end{eqnarray}
%
%
%
%
%
%% RETURN HERE -- INDEX PROBLEM!
%\begin{eqnarray} \label{error_equation_VC_1}
%E_{n} & = &  \prod_{\mu=0}^{n-1} Q_{n} E_0  \,+\,  \dt \, \sum_{\nu=0}^{n-1}  \prod_{\mu=1}^{n-\nu-1} Q_{n-\mu-1} \vtau_{\nu} \nonumber \\
%& = &  \prod_{\mu=0}^{n-1} Q_{n} E_0  \,+\,  
%\dt \, \sum_{\nu=0}^{n-2}  \prod_{\mu=1}^{n-\nu-1} 
% \left ( A + \dt B f(t_n) + (\mu+1) O(\dt^2) \right)  \vtau_{\nu}  +
% \dt \vtau_{n-1}
% \nonumber  
% \end{eqnarray}
 The first term is negligible because we assume that the initial error can be made arbitrarily small, 
 and the final term is clearly of order $\dt \vtau_{n-1}$.
 Using \eqref{EIS_Q_vc_def_4}, \eqref{EIS_Q_vc_def_5} and the same analysis as in \eqref{construction_criteria_1}--\eqref{construction_criteria_3} 
 we have 
 \begin{eqnarray*} 
\left\| \dt \, \sum_{\nu=0}^{n-2}  \left( \prod_{\mu=\nu+1}^{n-1} \tilde{Q}_{\mu}  \right)  \vtau_{\nu}  \right\| & = &
\left\| \dt \, \sum_{\nu=0}^{n-2}  \left( \prod_{\mu=\nu+2}^{n-1} \tilde{Q}_{\mu}  \right)   \left( \tilde{Q}_{\nu+1} \vtau_{\nu} \right) \right\| \\
& \leq & \dt \, \sum_{\nu=0}^{n-2}  \left\|  \prod_{\mu=\nu+2}^{n-1} \tilde{Q}_{\mu}  \right\|   \left\| \tilde{Q}_{\nu+1} \vtau_{\nu}  \right\| \\
& \leq & \dt \, \sum_{\nu=0}^{n-2} O\left( 1+O(\dt) \right)^{n- \nu -2}   O(\dt) \left\| \vtau_{\nu}  \right\| \\
& \leq & O(\dt) \max_{\nu=0, . . .,  n-2}   \left\| \vtau_{\nu}  \right\|  .\\
 \end{eqnarray*}
Putting these all together we have
 \begin{equation}\label{error_equation_VC_10}
 \left\|   E_n\right\| = O(\dt)  \max_{\nu=0, . . .,  n-1} \left\|   \vtau_{\nu}\right\| .
 \end{equation}

This simple proof shows that even for the variable coefficient case, the schemes constructed as described above have a higher order error
than would be expected from the truncation error. In the next subsection we extend this analysis to the general nonlinear case.

\subsection{Nonlinear equations} \label{section:EIS_NL_1}
Finally,  we analyze  the behavior of methods satisfying the assumptions in Section \ref{section:EIS_CC_1}
when applied  to  nonlinear problems. 
Consider the nonlinear equation
\begin{eqnarray}\label{EIS_nl_10}
&& u_t =  f(u(t), t)   \;\;,\;\;\;\;\;  t \ge 0 \nonumber \\
&& u(t=0) = u_0 \;
\end{eqnarray}
where $f(u, t)$ assumed to be analytic in $u$ and $t$. We now use the scheme
\begin{equation} \label{EIS_multistep_nl_1}
V_{n+1} = A V_n + \dt B \,
\left (  \begin{array}{cccccc}
f \left( v_{n+(s-1)/s},  t_{n+(s-1)/s} \right ) \\
\vdots \\
 f \left( v_n, t_{n} \right )
\end{array} \right )
\end{equation}
where the matrices $A$ and $B$ are as constructed above for the constant coefficients problem.

As defined in \eqref{truncation_error_1}, the exact solution to \eqref{EIS_nl_10} and the truncation error are related by
\begin{equation} \label{EIS_multistep_nl_2}
U_{n+1} = A U_n + \dt B \,
\left (  \begin{array}{cccccc}
f \left( u_{n+(s-1)/s},  t_{n+(s-1)/s} \right ) \\
\vdots \\
 f \left( u_n, t_{n} \right )
\end{array} \right ) + \dt \vtau_n.
\end{equation}
Note that  by Taylor expansion
\begin{equation} \label{EIS_multistep_nl_2.5} \nonumber
f \left( v_{\nu},  t_{\nu} \right ) =f \left( u_{\nu},  t_{\nu} \right ) +f _u\left( u_{\nu},  t_{\nu} \right ) (v_{\nu} -u_{\nu})+ r(v_{\nu} -u_{\nu})  \;,\;\
\end{equation}
where$f_u (u,t)= \partial f(u,t)/ \partial u$ and $|r(v_{\nu} -u_{\nu}) | \le c_1|v_{\nu} -u_{\nu}|^2$.
Subtracting \eqref{EIS_multistep_nl_1} from \eqref{EIS_multistep_nl_2} and  assuming that $E_n = U_n-V_n \ll 1$ gives
\begin{equation} \label{EIS_multistep_nl_3}
E_{n+1} = A E_n - \dt B \,
\left (  \begin{array}{cccccc}
f_u \left( u_{n+(s-1)/s},  t_{n+(s-1)/s} \right ) \\
& \ddots \\
&& \;\;\;\; f_u \left( u_n, t_{n} \right )
\end{array} \right ) E_n + \dt \vtau_n +\dt R (E_n) %O\left(\dt E_n^2 \right ) \;,
\end{equation}
where $\| R (E_n) \| \le c_1 \|E_n\|^2$. Equation \eqref{EIS_multistep_nl_3} means that as long as $O(E_n^2 ) \ll  O(\vtau_n)$,  
the equation for the error $E_n$ can be analyzed in essentially the same way as for the
 linear  variable coefficient case, and the same estimates hold.

In order to evaluate the time interval in which $O(E_n^2 ) \ll  O(\vtau_n)$ we note that although the term $R (E_n)$  
in \eqref{EIS_multistep_nl_3} is not a non-homogeneous term but rather a function  of $E_n$, 
we can still use the approach used in  \cite[Theorem 5.1.2]{gustafsson1995time})
to prove stability for a perturbed solution operator. As in \cite[Theorem 5.1.2]{gustafsson1995time}), we use the
discrete version of Duhamel's principle to obtain
\begin{eqnarray} \label{EIS_multistep_nl_10}
E_{n} & = &  \prod_{\mu=0}^{n-1} \hat{Q}_{\mu}  E_0  \,+\,  \dt \, \sum_{\nu=0}^{n-1}  
\left( \prod_{\mu=\nu+1}^{n-1} \hat{Q}_{\mu}  \right)  \vtau_{\nu}  + \dt \sum_{\nu=0}^{n-1}  
\left( \prod_{\mu=\nu+1}^{n-1} \hat{Q}_{\mu}  \right)   R (E_{\nu}) \nonumber \\
%& = &  \prod_{\mu=0}^{n-1} Q_{\mu} E_0  \,+\,  
%\dt \, \sum_{\nu=0}^{n-2}  \prod_{\mu=\nu+1}^{n-1}
% \left ( \tilde{Q}_{\mu} +  O(\dt^2) \right)  \vtau_{\nu}  +
% \dt \vtau_{n-1}
% \nonumber  
 \end{eqnarray}
where
\begin{equation} \label{EIS_multistep_nl_20}
\hat{Q}_{n} = A  - \dt B \,
\left (  \begin{array}{cccccc}
f_u \left( u_{n+(s-1)/s},  t_{n+(s-1)/s} \right ) \\
& \ddots \\
&& \;\;\;\; f_u \left( u_n, t_{n} \right )
\end{array} \right )\;.
\end{equation}
Taking the norm of \eqref{EIS_multistep_nl_10} and using the triangle inequality we obtain
\begin{eqnarray} \label{EIS_multistep_nl_30}
\left \| E_{n} \right \|& \le &  \left \|  \prod_{\mu=0}^{n-1} \hat{Q}_{\mu}  E_0 \right \| \,+\,  \left \|  \dt \, \sum_{\nu=0}^{n-1}  
\left( \prod_{\mu=\nu+1}^{n-1} \hat{Q}_{\mu}  \right)  \vtau_{\nu}  \right \|+ \left \| \dt \sum_{\nu=0}^{n-1}  
\left( \prod_{\mu=\nu+1}^{n-1} \hat{Q}_{\mu}  \right)   R (E_{\nu})  \right \| .\nonumber \\
 \end{eqnarray}

As in the linear case we assume that the initial error, $E_0$ is arbitrary small, so the first term is negligible.
If $\hat{Q}_{\nu+1}$ is constructed such that $\| \hat{Q}_{\nu+1} \vtau_{\nu}  \| = \dt O(\vtau_{\nu})$ then using 
the same analysis as in variable coefficient case the second term in \eqref{EIS_multistep_nl_30} is less 
or equal to  $\dt c_0 \phi_h^*(c, \, t_n) \max_{\nu=0, . . .,  n-1} \left\|   \vtau_{\nu}\right\| $.  
As to the third term, the same arguments can be used to show that it is bounded by
\begin{equation} \label{EIS_multistep_nl_60}
\left \| \dt \sum_{\nu=0}^{n-1}  
\left( \prod_{\mu=\nu+1}^{n-1} \hat{Q}_{\mu}  \right)   R (E_{\nu})  \right \|  \le 
c_1 \phi_h^*(c, \, t_n)  \left  \| E_{n}  \right \| ^2 \;,
\end{equation}
so that \eqref{EIS_multistep_nl_30} (with the substitution of \eqref{EIS_multistep_nl_60} for the final term)
can be re-arranged to obtain
 \begin{eqnarray} \label{EIS_multistep_nl_70}
\left \| E_{n} \right \|  \left( 1-  c_1 \phi_h^*(c, \, t_n)  \left  \| E_{n}  \right \| \right)& \le &  \dt c_0 \phi_h^*(c, \, t_n) 
\max_{\nu=0, . . .,  n-1} \left\|   \vtau_{\nu}\right\| .  %\nonumber \\
 \end{eqnarray}
 If $ c_1 \phi_h^*(c, \, t_n)  \left  \| E_{n}  \right \|<1/2$, we obtain
  \begin{eqnarray} \label{EIS_multistep_nl_80}
\left \| E_{n} \right \|  & \le &  2 \dt c_0 \phi_h^*(c, \, t_n) \max_{\nu=0, . . .,  n-1} \left\|   \vtau_{\nu}\right\|   %\nonumber \\
 \end{eqnarray}
 This estimate holds as long as
   \begin{eqnarray} \label{EIS_multistep_nl_90}
c_1 \phi_h^*(c, \, t_n) \left \| E_{n} \right \|  & \le &  2 \dt c_0 c_1\left( \phi_h^*(c, \, t_n) \right)^2  \max_{\nu=0, . . .,  n-1} 
\left\|   \vtau_{\nu}\right\|  \le \frac{1}{2}, %\nonumber \\
 \end{eqnarray}
 which is satisfied for all times $t_n$ such that   $\dt \phi_h^*(c, \, t_n)= O(1)$.
 
 Therefore
  \begin{equation}\nonumber \label{EIS_multistep_nl_100}
 \left\|   E_n\right\| = O(\dt)  \max_{\nu=0, . . .,  n-1} \left\|   \vtau_{\nu}\right\| .
 \end{equation}
for the nonlinear case as well. 
 
%
%We justify the assumption $E_n = U_n-V_n \ll 1$ recursively:
%we know that $E_0 \ll 1$, and the relation   \eqref{EIS_multistep_nl_3} above shows that whenever 
%$E_n \ll 1$ we also have $E_{n+1} \ll 1$. \\
%
%WE NEED TO GO OVER IT ADI.
%%S_h\left( t_n, t_{\nu}\right ) = \prod_{\mu=1}^{n-\mu} Q_{n-\mu} \;,\;\;\;S_h\left( t_n, t_n\right )  = I

%%% 
% u_{n+(s-1)/s}, \ldots,  u_{n+1/s},  u_n\right

%We assume that $Q$ can be diagonalzed, i.e. that there is a  that it has $s$ linearly independent eigenvectors for all values of $0\le k<k_0$. Note that in this case  $S_h\left( t_n, t_{\nu}\right ) = Q^{n-\nu}$. The stability of the scheme implies that all the eigenvalues of $Q$, $z_j$, satisfies  $ | z_j |\le 1+O(k)$.

%Since $u(t_{n+1}) \, u(t_n)+O(k)$ then if the numerical solution $V_n$ converges to the analytic solution $U_n$ then one of the eigenvalues of $Q$ should be equal to $1+O(k)$ and one eigenvector should have the form:
%\begin{equation} \label{2.20}
%\left( 1+O(k), \ldots, 1+O(k)\right)^T
%\end{equation}

%%%%%%%%%%%%%%%%%%%%%%%%%%%%%%%%%%%%%%%%%%%

\section{Some Error Inhibiting Explicit Schemes}  \label{section:examples}

In the previous section we define sufficient conditions for methods of the form
\begin{equation} 
V_{n+1} = Q V_n
\end{equation}
where 
\[ Q = A + \dt B f \] 
to be error inhibiting. These are
\begin{enumerate}
\item[\bf C1.] ${\bf rank} (A)=1$.
\item[\bf C2.]  Its non-zero eigenvalue is equal to $1$ and its corresponding eigenvector is 
\[ \left( 1, \ldots, 1\right)^T.\]
\item[\bf C3.] $A$ can be diagonalized. 
\item[\bf C4.] The matrices $A$ and $B$ are constructed such that when the  local truncation error  is multiplied by the discrete solution operator  
we have the bound:
\[ \left \|  Q \vtau_{\nu} \right \|  \leq O(\dt)  \left \| \vtau_{\nu} \right \| . \]
This is accomplished by requiring the local truncation error to live in the space of the eigenvectors of $A$ that correspond to the 
zero eigenvalues.
\end{enumerate}
In this section we present several schemes which were constructed using the approach presented in the previous section. 
In Section \ref{s2p3method}, we present a block one-step method that evolves two steps ($v_n$ and $v_{n+\frac{1}{2}}$)
 to obtain the next two steps ($v_{n+1}$ and $v_{n+\frac{3}{2}}$). 
This method has truncation error \eqref{truncation_error_1} that is second order, while its global order \eqref{error} is third order.  
We demonstrate that the expected convergence  rate is attained on several  sample nonlinear problems. In this section we also show that
a typical Type 3 DIMSIM method (derived in \cite{butcher1993a}) that satisfies the first three conditions above but not the fourth,
has truncation error of order two, and its global error is of the same order. This demonstrates the importance of condition {\bf C4}.

Next, in Section \ref{s3p4method} we present a block one-step  method that evolves three 
steps $v_n$,  $v_{n+\frac{1}{3}}$ and$v_{n+\frac{2}{3}}$ 
to obtain $v_{n+1}$, $v_{n+\frac{4}{3}}$ and $v_{n+\frac{5}{3}}$. 
This method has truncation error \eqref{truncation_error_1} that is third order, while its global order  \eqref{error}  is fourth order, as we 
demonstrate on several  sample problems.
Finally, to show that the methods in each class are not unique, 
we present two other methods of this type and show that their global error is of one order higher than the local truncation error
on a sample nonlinear system.

\subsection{A third order error inhibiting  method with $s=2$.} \label{s2p3method}
In this subsection we define an explicit block one-step with $s=2$ 
that satisfies the conditions {\bf C1 --  C4} above. This method
 takes the values of the solution at the times $t_n$ and
$t_{n+\frac{1}{2}}$ and obtains the solution at the time-level $t_{n+1}$ and $t_{n+\frac{3}{2}}$.
 The  exact solution vector for this problem is
\[ U_n = \left (  u(t_{n+1/2}),   u(t_n)  \right )^T\]
and, similarly, the corresponding vector of numerical approximations is
\[  V_n = \left( v_{n+1/2},   v_n\right )^T. \]
The scheme is given by:
\begin{equation} \label{EIS_2_step_10}
V_{n+1}  \,=  \, \frac{1}{6}\left( \begin{array}{ccc}
-1& \;\; & 7 \\
-1&  & 7
\end{array} \right ) V_n + 
\frac{\dt}{24}\left( \begin{array}{ccc}
55 & \;&-17 \\
25 & & 1
\end{array} \right ) 
 \, \left (  \begin{array}{cccccc}
f \left( v_{n+1/2},  t_{n+1/2} \right ) \\
 f \left( v_n, t_{n} \right )
\end{array} \right ),
\end{equation}
and  has truncation error 
\begin{equation} \label{EIS_2_step_20}
\vtau_{n}  \,=  \, \frac{23}{576}\left( \begin{array}{ccc}
7 \\
 1
\end{array} \right ) \frac{d^3} {dt^3}  u(t_n) \,\Delta t^2 \, + \, O(\Delta t^3) \;.
\end{equation}

The matrix $A$ can be diagonalized as follows:
\begin{equation} \label{EIS_2_step_30}
A  \,=  \, \frac{1}{6}\left( \begin{array}{ccc}
-1& \;\; & 7 \\
-1&  & 7
\end{array} \right ) \,=  \,  \frac{1}{6} \left( \begin{array}{ccc}
1& \;\; & 7 \\
1&  & 1
\end{array} \right )  \left( \begin{array}{ccc}
1& \;\; &  \\
&  & 0
\end{array} \right ) \left( \begin{array}{ccc}
-1& \;\; & 7 \\
1&  & -1
\end{array} \right )  \;.
\end{equation}
Observe that the leading order of the truncation error \eqref{EIS_2_step_20} is in the space of the second eigenvector of $A$, the one that corresponds
to the zero eigenvalue. Also, as was pointed out in Remark \ref{remark:truncation_error_rem},  $\vtau_n$ depends only on this eigenvector of $A$ and a multiple that
is not directly dependent on $f$ but only on the third derivative of the solution $u$.
This underscores the analysis in Sections  \ref{section:EIS_VC_1} and \ref{section:EIS_NL_1}
that demonstrates that the error inhibiting property carries through for variable coefficient and nonlinear 
problems.

To study the behavior of the global error 
%we observe the leading order, $O(\Delta t^0)$, the eigenvalue of $Q$
%which corresponds to the exact solution is one,  while the eigenvalue corresponding to the truncation error is zero. 
%For the next order, $O(\Delta t^1)$, 
we use the fact shown in Section \ref{section:EIS_NL_1} that even for a nonlinear equation 
 it is sufficient to analyze the matrix 
\begin{equation} \label{EIS_2_step_40}
Q \,=  \, A + \dt \,B\,f
\end{equation}
where $f$ is a constant. 
In this case:
\begin{eqnarray}\label{EIS_2_step_50}
Q   &= &  \frac{1}{6} \left( \begin{array}{ccc}
1 + \frac{{f} \dt}{2} + \frac{{f}^2 \dt^2}{8}+O(\dt^3) & \;\; & 7 + 36 {f}  \dt + 228 {f}^2 \dt^2 +O(\dt^3) \\
1&  & 1
\end{array} \right )  \nonumber \\
&& \left( \begin{array}{ccc}
1 + {f} \dt + \frac{{f}^2 \dt^2}{2} + \frac{{f}^3 \dt^3}{6}+O(\dt^4)& \;\; &  \\
&  & \frac{4 {f} \dt}{3} - \frac{{f}^2 \dt^2}{2} - \frac{{f}^3 \dt^3}{6}+O(\dt^4)
\end{array} \right ) \nonumber \\
&&\left( \begin{array}{ccc}
-1 + \frac{71 {f} \dt}{12} + \frac{107 {f}^2 \dt^2}{36} +O(\dt^3)& \;\; & 7 - \frac{65 {f} k}{12} - \frac{209 {f}^2 \dt^2}{36}+O(\dt^3) \\
1 - \frac{71 {f} \dt}{12} - \frac{107 {f}^2 \dt^2}{36}+O(\dt^3)&  & -1 + \frac{65 {f} k}{12} + \frac{209 {f}^2 \dt^2}{36}+O(\dt^3)
\end{array} \right )  \nonumber \\
\end{eqnarray} 

Recall that, neglecting the initial error $E_0$, we can say that the global error is  \eqref{error}
\[ E_n = \dt \sum_{\nu=0}^{n-1} Q^{n-\nu-1} \tau_\nu \]
Putting together equations \eqref{EIS_2_step_20} and \eqref{EIS_2_step_50} we see that each term
 $ Q \,\vtau_\nu$ contributes to the error in two ways:
 \begin{itemize}
\item The first contribution is due to the fact that $ \vtau_\nu$ is almost co-linear with the second eigenvector $\psi_2$. 
The order of this  contribution is 
\[ |z_2|  \| \psi_2 \vtau_\nu \| = O(\dt) \cdot O(\| \dt \vtau_\nu \|) =  O(\dt^3) \]
 where the term $|z_2|$ is the second eigenvalue which is of order $O(\dt)$. 
\item The second contribution to the error comes from the component of  $ \vtau_\nu$ 
that is a multiple $\gamma_1$ of the first eigenvector $\psi_1$, 
\[ |z_1| \| \gamma_1 \psi_1  \vtau_\nu \|  = O(\dt) \cdot O(\| \vtau_\nu \|) = O(\dt^3) \]
the term $\gamma_1$ is of $O(\dt)$ because $ \vtau_\nu$ lives mostly in the space of $\psi_2$.
\end{itemize} 
While each of the terms in  $ \dt Q \,\vtau_\nu$ has  order $O(\dt^2) \cdot O(\| \vtau_\nu \|) = O(\dt^4)$,
as the method is evolved forward, the errors accumulate over time,  and
sum of all contributions from all the times  gives us a global error of order  $O(\dt) \cdot O(\| \vtau_n \|) = O(\dt^3)$.

\begin{figure*}[t!]
    \centering
    \begin{subfigure}[t]{0.5\textwidth}
        \centering
        \includegraphics[width=0.925\textwidth]{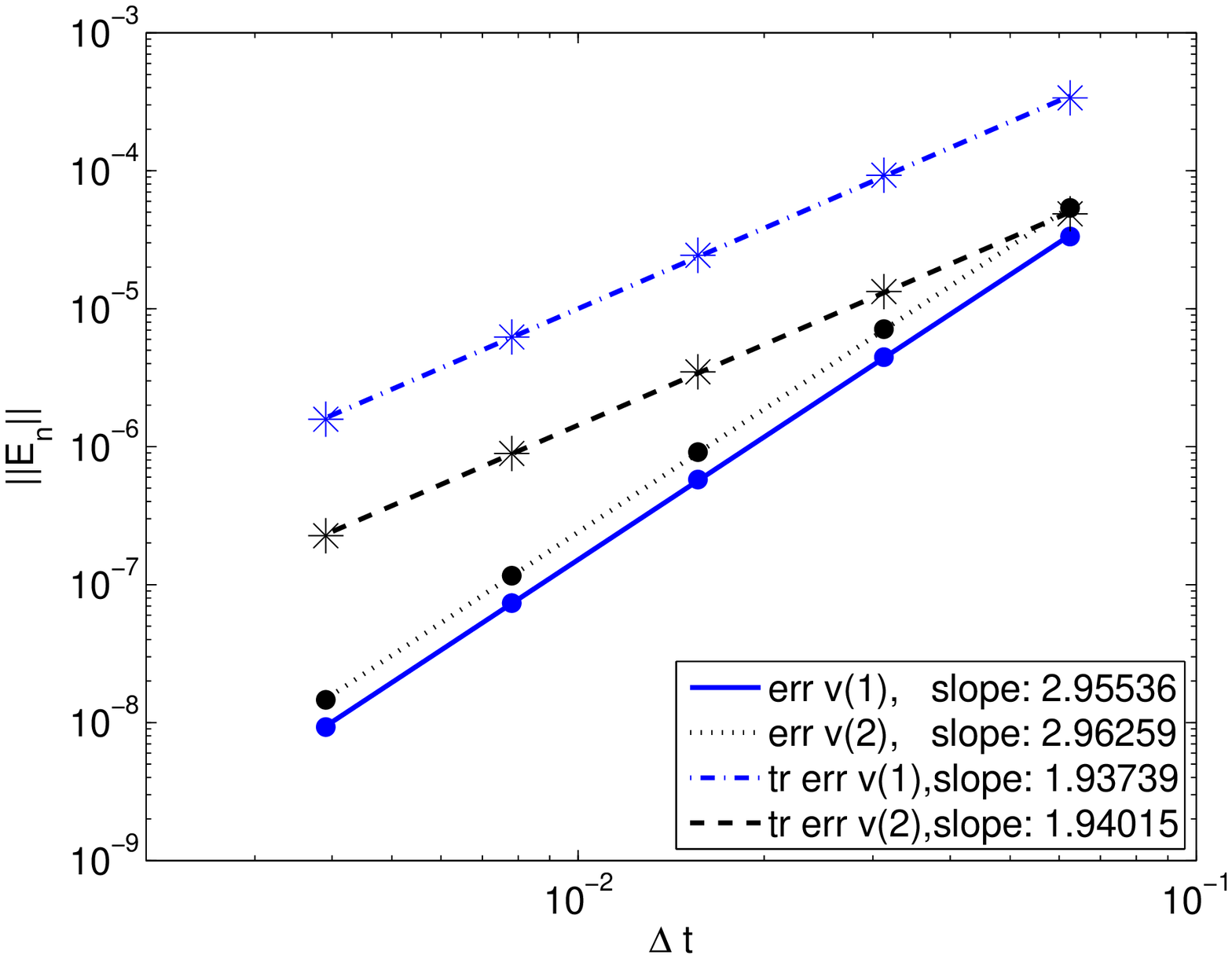}
        \caption{}
    \end{subfigure}%
    ~ 
    \begin{subfigure}[t]{0.5\textwidth}
        \centering
      \includegraphics[width=\textwidth]{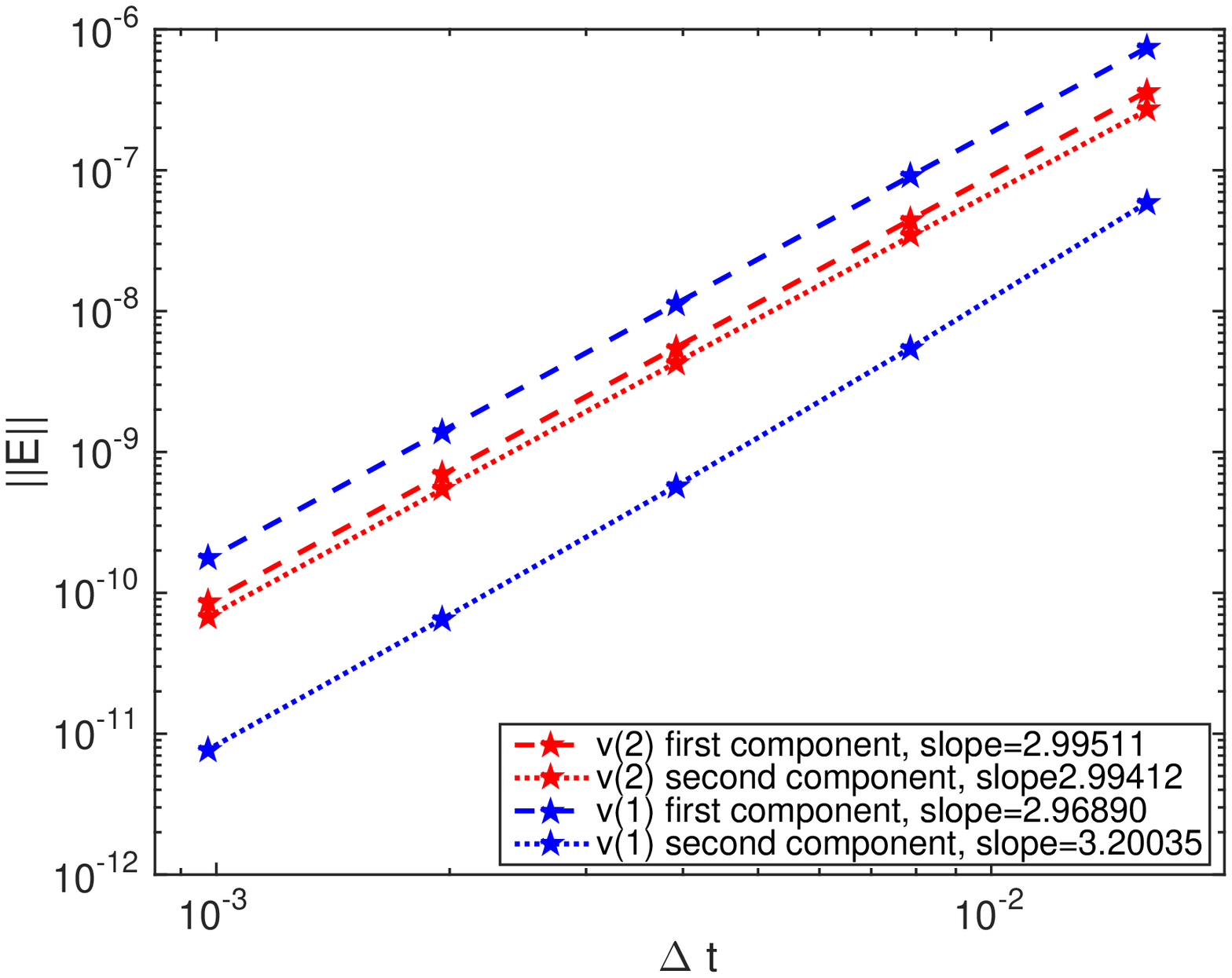}
        \caption{}
    \end{subfigure}
    \caption{Convergence plots using the scheme \eqref{EIS_2_step_10}. (a) The errors and truncation errors
vs. $\Delta t$, for several values of $\Delta t$, for the numerical solution of \eqref{EIS_2_step_example_1_10}.
(b) The errors  vs. $ \dt$ for each component of the solution,
computed for several values of $\dt$, for the numerical solution of the van der Pol equation \eqref{VDP_1}.}
 \label{fig:EIS_2_step_example}
\end{figure*}

{\bf Example 1a}: To demonstrate that this method indeed performs as designed we study its behavior
 on a  nonlinear scalar equation of the form:
\begin{eqnarray}\label{EIS_2_step_example_1_10}
& & u_t=   -u^2 \; = \; f(u)   \;\;,\;\;\;\;\;  t \ge 0 \nonumber \\
& & u(t=0) = 1 \;.
\end{eqnarray}
We evolve the solution of this equation to time $T=1$ using the scheme \eqref{EIS_2_step_10}. 
The initial steps are computed exactly. %The truncation error is computed for comparison, and 
The plots of the errors and the truncation errors are presented in Figure \ref{fig:EIS_2_step_example}(a). 
Both errors are shown for the first component, $v_{n+1/2}$ (denoted v(1) in the legend) and the second component,
$v_{n}$ (denoted v(2) in the legend).
Clearly, although the truncation error is only second order (denoted tr  err v(1) and tr  err v(2) in the legend), the global error is third order, as predicted by the theory.

%\begin{figure}[h]
%\subfloat[]{}
%\subfloat[]{\includegraphics[width=0.5\textwidth]{./Figures/vdp_2step.eps}}
%\caption{Convergence plots using the scheme \eqref{EIS_2_step_10}. On the left are the errors and truncation errors
%in the  numerical solution of \eqref{EIS_2_step_example_1_10} vs. $\Delta t$, for several values of $\Delta t$.
%On the right are the errors  vs. $ \dt$ for each component of the solution,
%computed for several values of $\dt$, for the numerical solution of the van der Pol equation \eqref{VDP_1}.
%}
%
%\end{figure}

{\bf Example 1b}: It is important that the method will perform as designed on a nonlinear system as well. To demonstrate this,
we solve the the van der Pol system 
\vspace{-.1in}
\begin{eqnarray}\label{VDP_1}
u^{(1)}_t & = & u^{(2)} \nonumber \\
u^{(2)}_t & = & 0.1  [1-  (u^{(1)})^2 ] u^{(2)}-  u^{(1)} 
\end{eqnarray}
using the same scheme \eqref{EIS_2_step_10}.  As this is a system, it is important that both components are examined.
Thus, the vector of the numerical solution has two components for the time level $t_n$, denoted by v(2), and two components
for the time level $t_{n+\frac{1}{2}}$, denoted by v(1).
In Figure \ref{fig:EIS_2_step_example}(b)  the convergence plot of the  components of $u^{(1)}$ and $u^{(2)}$ are presented. 
Once again, we see that the convergence rate is indeed third order.

%\begin{figure}[h]
%\begin{center}
%\begin{tabular}{ccccc}
%\includegraphics[width=0.7\textwidth]{./Figures/vdp_2step.eps}&
%\end{tabular}
%\end{center}
%\caption{Convergence plots for the numerical solution of the van der Pol equation using using scheme \eqref{EIS_2_step_10}. The plots are for $\log_{10} \|\vE\|$ vs. $\log_{10} \dt$ for each component of the solution,
%computed for several values of $\dt$}
% \label{fig:EIS_2_step_example_1_11}
%\end{figure}

%{\bf Remark:} {\em It is important to note that not all Type 3 DIMSIM methods  have the EIS property! The property that the local truncation error
%lives in the space spanned by the eigenvectors of $A$ that correspond to the zero eigenvalues is needed for the error inhibiting behavior
%to occur, and this property  is not generally satisfied. To observe this, we study the DIMSIM scheme of types 3 presented by   J. C. Butcher in \cite{butcher1993a}. }

\begin{remark}\label{remark:dimsim_vs_EIS} 
 It is important to note that not all Type 3 DIMSIM methods  have the EIS property! The property that the local truncation error
lives in the space spanned by the eigenvectors of $A$ that correspond to the zero eigenvalues is needed for the error inhibiting behavior
to occur, and this property  is not generally satisfied. To observe this, we study the DIMSIM scheme of types 3 presented by   J. C. Butcher in \cite{butcher1993a}.
\end{remark}

Consider the  scheme 
\begin{equation} \label{EIS_2_step_Butcher_10}
  \left( \begin{array}{ccc}
v_{n+2}\\
v_{n+1}
\end{array} \right )  \,=  \, \frac{1}{4}\left( \begin{array}{ccc}
7& \;\; & -3 \\
7&  & -3
\end{array} \right ) \left( \begin{array}{ccc}
v_{n+1}\\
v_{n}
\end{array} \right )+ 
\frac{\dt}{8}\left( \begin{array}{ccc}
9 & \;&-7 \\
-3 & & -3
\end{array} \right ) 
 \, \left (  \begin{array}{cccccc}
f \left( v_{n+1},  t_{n+1} \right ) \\
 f \left( v_n, t_{n} \right )
\end{array} \right )
\end{equation}
given in \cite{butcher1993a}.
This scheme has truncation error 
\begin{equation} \label{EIS_2_step_Butcher_20}
\vtau_{n}  \,=  \, \frac{1}{48}\left( \begin{array}{ccc}
23 \\
 3
\end{array} \right ) \frac{d^3} {dt^3}  u(t_n) \,\Delta t^2 \, + \, O(\Delta t^3) \;.
\end{equation}

The matrix $A$ can be diagonalized as follows:
\begin{equation} \label{EIS_2_step_Butcher_30}
A  \,=  \, \frac{1}{4}\left( \begin{array}{ccc}
7& \;\; & -3 \\
7&  & -3
\end{array} \right )\,=  \,  \left( \begin{array}{ccc}
1& \;\; & 3/7 \\
1&  & 1
\end{array} \right )  \left( \begin{array}{ccc}
1& \;\; &  \\
&  & 0
\end{array} \right )  \frac{1}{4}\left( \begin{array}{ccc}
7& \;\; & -3 \\
-7&  & 7
\end{array} \right )  \;.
\end{equation}
The  truncation error $\vtau_{n}$ can be written as a linear combination of the two eigenvectors of $A$ as follows:
\begin{equation} \label{EIS_2_step_Butcher_40}
\vtau_{n}  \,=  \, \left[ \frac{19}{24}\left( \begin{array}{ccc}
1 \\
 1
\end{array} \right ) - \frac{35}{48}\left( \begin{array}{ccc}
3/7 \\
 1
\end{array} \right )\right ] \frac{d^3} {dt^3}  u(t_n) \,\Delta t^2 \, + \, O(\Delta t^3) \;.
\end{equation}
Unlike the EIS scheme \eqref{EIS_2_step_10}, here the first term in this expansion is of the order of $O(\vtau_{n})= O(\Delta t^2)$. 
Therefore a term of the order of $\Delta t O(\vtau_{n})= O(\Delta t^3)$ is accumulated at each time step,
so that the global error is second order.

\begin{figure*}[t!]
    \centering
    \begin{subfigure}[t]{0.5\textwidth}
        \centering
        \includegraphics[width=0.95\textwidth]{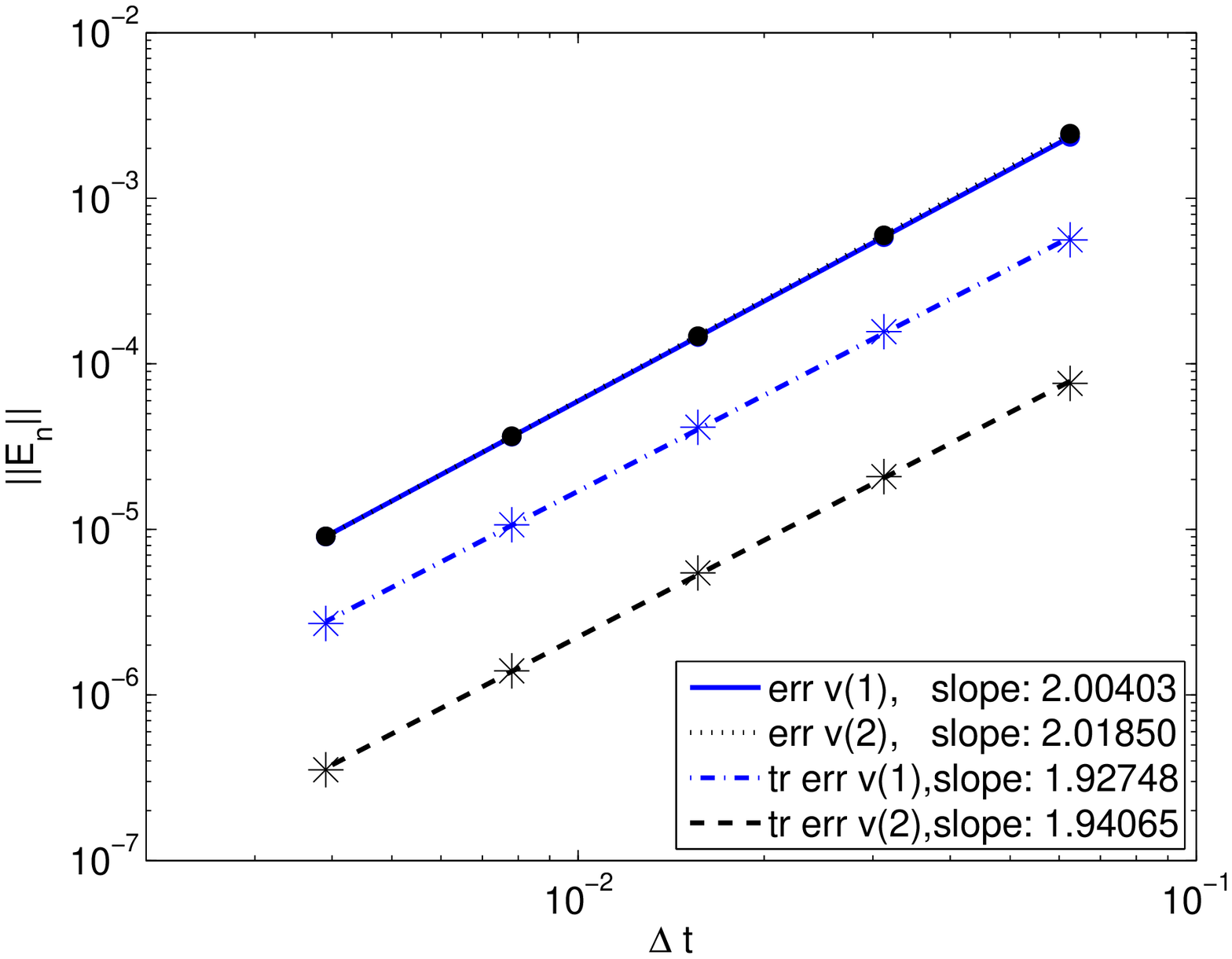}
        \caption{}
    \end{subfigure}%
    ~ 
    \begin{subfigure}[t]{0.5\textwidth}
        \centering
      \includegraphics[width=1.025\textwidth]{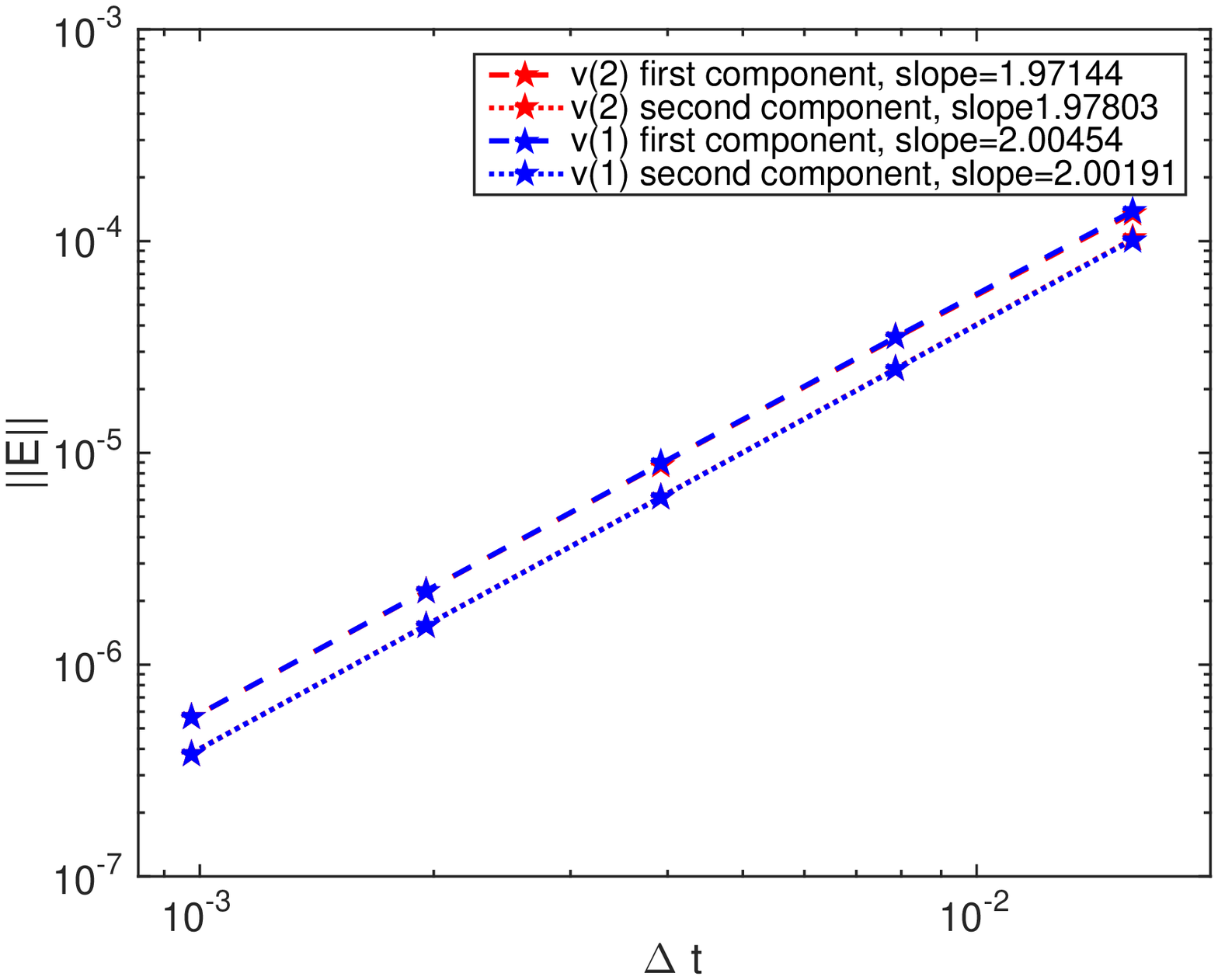}
        \caption{}
    \end{subfigure}
    \caption{Convergence plots using Butcher's scheme \eqref{EIS_2_step_Butcher_10}. 
    (a) The errors and truncation errors
vs. $\Delta t$, for several values of $\Delta t$, for the numerical solution of \eqref{EIS_2_step_example_1_10}.
Note that the errors for v(1) and  v(2) are virtually identical so these error lines coincide.
(b) The errors  vs. $ \dt$ for each component of the solution,
computed for several values of $\dt$, for the numerical solution of the van der Pol equation \eqref{VDP_1}.
Note that for this problem as well the behavior of this method on both components is virtually identical, so the error lines  for 
each component of the solution coincide.
Both the local truncation errors and the global errors are second order: this is not an error inhibiting scheme.
    }
 \label{fig:EIS_2_step_Butcher_example}
\end{figure*}

We note that both this method \eqref{EIS_2_step_Butcher_10} and our error inhibiting method 
 \eqref{EIS_2_step_10} satisfy the order conditions in Theorem 3.1 of \cite{butcher1993a} only up to second
 order ($p=2$). However, as we see in Figure \ref{fig:EIS_2_step_Butcher_example}, 
 when the method \eqref{EIS_2_step_Butcher_10} is used to simulate the solution 
 of the problems \eqref{EIS_2_step_example_1_10} and \eqref{VDP_1} we have second order 
 accuracy,  while the error inhibiting method  \eqref{EIS_2_step_10} gave third order accuracy 
 (Figure \ref{fig:EIS_2_step_example}).

\subsection{A fourth order error inhibiting method with $s=3$.} \label{s3p4method}
In this subsection we present an error inhibiting method with $s=3$ that takes the values of the solution at the 
times $t_n$, $t_{n+\frac{1}{3}}$, and  $t_{n+\frac{2}{3}}$ and uses these three values to obtain the solution 
at the time-level $t_{n+1}$,  $t_{n+\frac{4}{3}}$, and $t_{n+\frac{5}{3}}$. The exact solution vector 
is given by
\[ U_n = \left (  u(t_{n+2/3}), u(t_{n+1/3}),   u(t_n) \right )^T, \]
and the corresponding vector of numerical approximations is 
\[ V_n = \left( v_{n+2/3}, v_{n+1/3},   v_n\right )^T.\]
Consider the error inhibiting scheme

\begin{eqnarray} \label{EIS_3_step_10}
V_{n+1}  &= &  \frac{1}{768}\left(
\begin{array}{ccccc}
467 & -1996 & 2297 \\
467 & -1996 & 2297 \\
467 & -1996  & 2297 \\
\end{array}
\right) V_n + \nonumber \\
&& \hspace{1cm}\frac{\dt}{1152}\left(
\begin{array}{ccc}
5439 &  -6046 & 3058 \\
2399 & -1694  & 1362 \\
703   &   354 & 626 \\
\end{array}
\right)
 \, \left (  \begin{array}{cccccc}
 f \left( v_{n+2/3},  t_{n+2/3} \right ) \\
f \left( v_{n+1/3},  t_{n+1/3} \right ) \\
 f \left( v_n, t_{n} \right )
\end{array} \right ),
\end{eqnarray}
which has a local truncation error of  third order,
\begin{eqnarray} \label{EIS_3_step_20}
\vtau_{n}  &= & \frac{1}{373248}\left( \begin{array}{ccc}
43699 \\
12787 \\
2227
\end{array} \right ) \frac{d^4} {dt^4}  u(t_n) \,\dt^3 \, + \, O(\dt^4) \nonumber \\
\nonumber \\
& & \approx \; \left( \begin{array}{ccc}
0.117078\\
 0.0342587\\
  0.00596654 \\
\end{array} \right ) \frac{d^4} {dt^4}  u(t_n) \,\dt^3 \, + \, O(\dt^4) \;.
\end{eqnarray}
However, it can be verified that for the linear case,  the product    
\[Q_n \vtau_{n}  = O(\dt \vtau_n) =  O(\dt^4) \; . \]
Given the analysis in  Section \ref{section:EIS_NL_1}
above, this result will carry over to the nonlinear case,
and  thus this method will have a  fourth order global error, despite the third order truncation error. 

\begin{figure*}[t!]
    \centering
    \begin{subfigure}[t]{0.5\textwidth}
        \centering
        \includegraphics[width=0.95\textwidth]{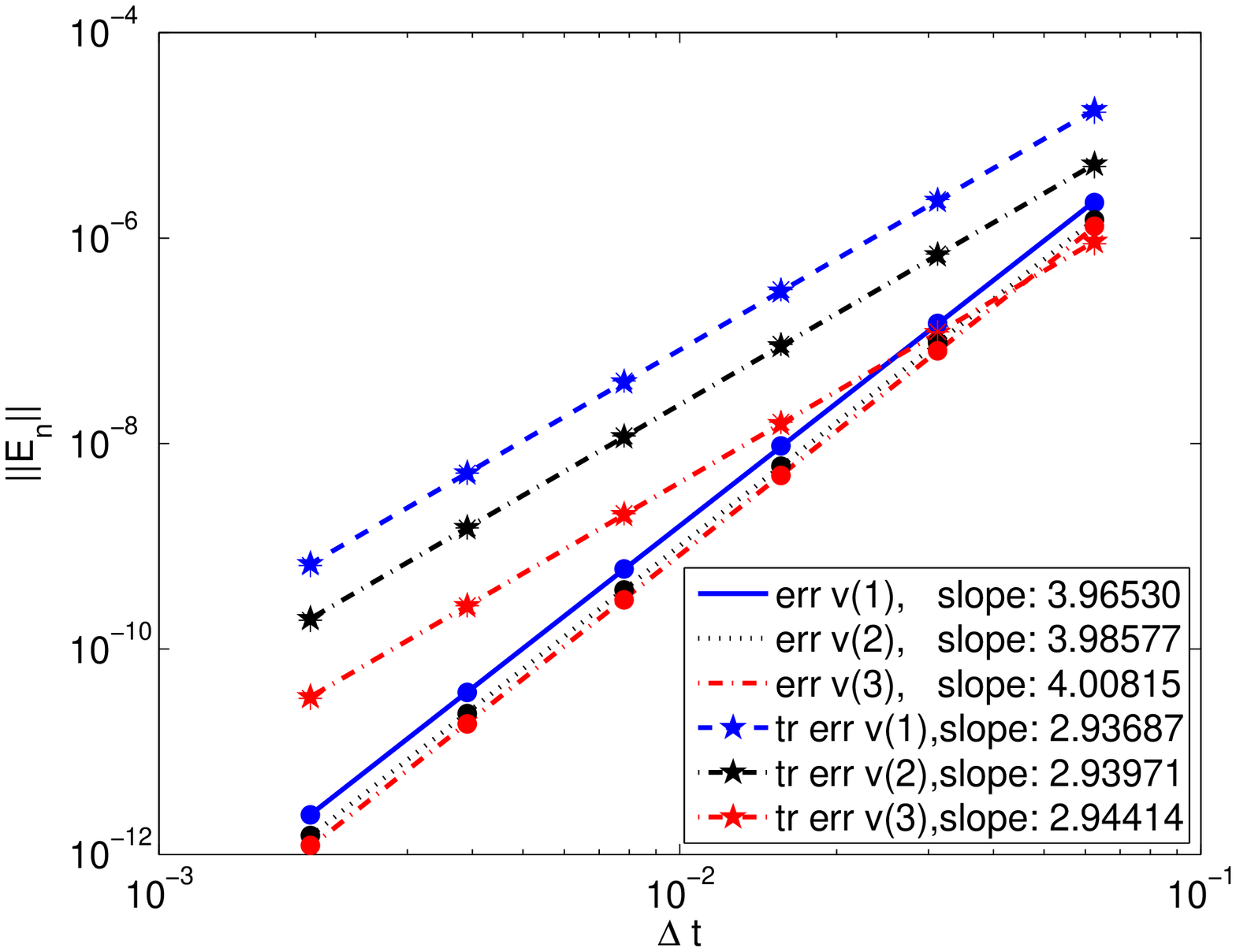}
        \caption{}
    \end{subfigure}%
    ~ 
    \begin{subfigure}[t]{0.5\textwidth}
        \centering
      \includegraphics[width=1.02\textwidth]{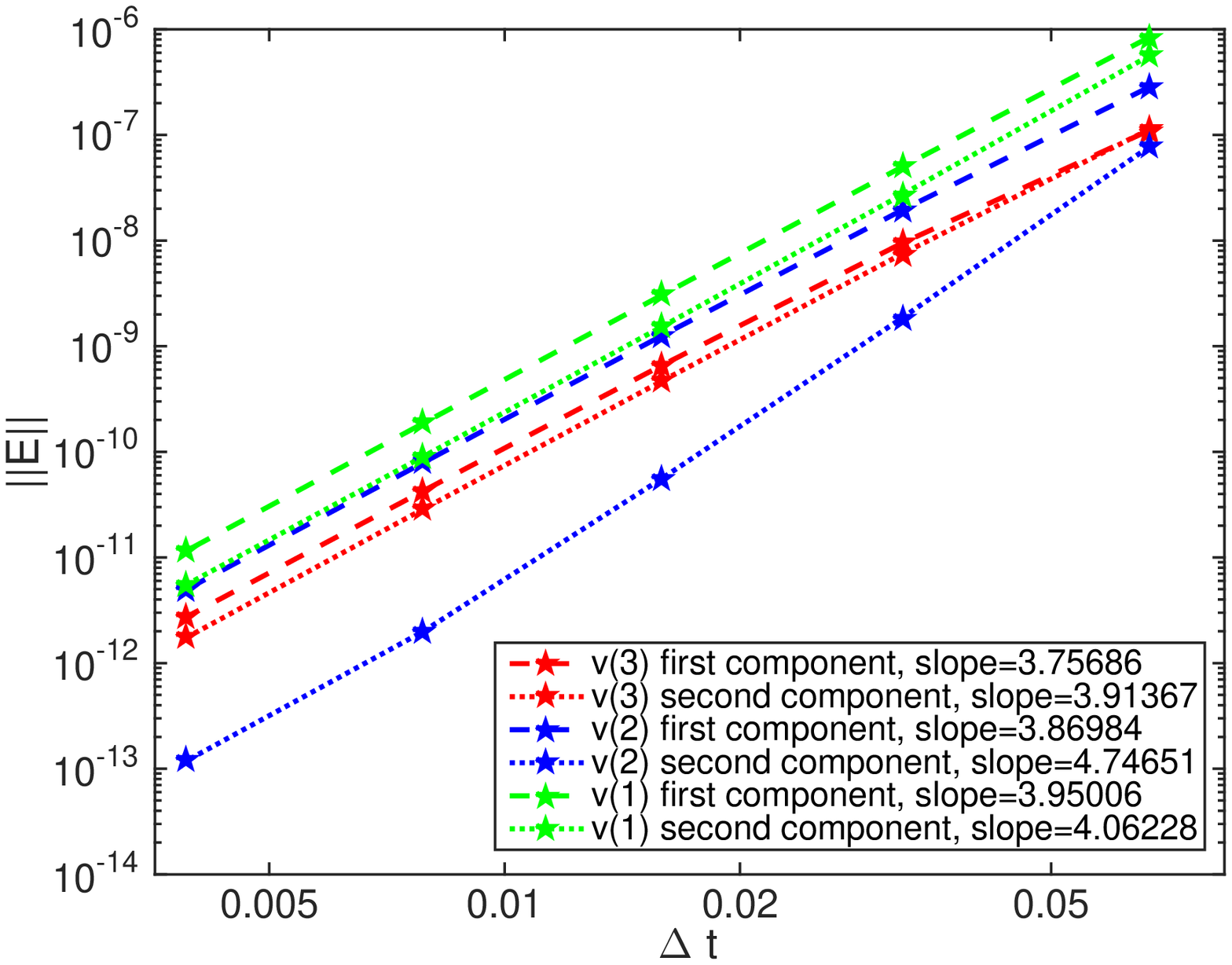}
        \caption{}
    \end{subfigure}
    \caption{Convergence plots using the scheme \eqref{EIS_3_step_10}. 
(a) The errors and truncation errors
vs. $\Delta t$, for several values of $\Delta t$, for the numerical solution of \eqref{EIS_2_step_example_1_10}.
(b) The errors  vs. $ \dt$ for each component of the solution,
computed for several values of $\dt$, for the numerical solution of the van der Pol equation \eqref{VDP_1}.
As expected, we observe fourth
order accuracy for the errors, although the truncation errors are third order. }
 \label{fig:EIS_3_step_example}
\end{figure*}

To demonstrate this result we revisit the two examples  \eqref{EIS_2_step_example_1_10} and \eqref{VDP_1} 
 in the previous subsection and  use the scheme \eqref{EIS_3_step_10} to evolve them forward in time. 
 The results, shown  in Figure \ref{fig:EIS_3_step_example}, are exactly as we expect: although the truncation errors 
(seen for the problem  \eqref{EIS_2_step_example_1_10} in Figure \ref{fig:EIS_3_step_example}(a))
are only third order, the errors are fourth order for both problems  \eqref{EIS_2_step_example_1_10} and
the van der Pol problem \eqref{VDP_1}. 

%We note that for the first problem, the errors for the third component feature an 
%unusual behavior for the points selected. The dip that can be seen is an interesting feature that does not
%detract from the accuracy of the method: in fact, the errors fall below the order of accuracy line one would expect.

\subsubsection{Other fourth order error inhibiting methods with $s=3$.}
%The methods above are not unique, in fact it is easy to obtain other methods using this approach. 
The methods above are not unique, in fact other methods can be derived using this approach. 
In this section we present two additional error inhibiting methods with $s=3$ that have local truncation error
that is third order  but demonstrate fourth order global error on a nonlinear system. 

The first  method is
\begin{eqnarray} \label{EIS_3_step_10b}
V_{n+1}  &= &  \frac{1}{1020}\left(
\begin{array}{ccccc}
 449 & -1966 &\;& 2537 \\
 449 & -1966 &\;& 2537 \\
 449 & -1966 &\;& 2537 \\
\end{array}
\right) V_n + \nonumber \\
&& \hspace{1cm}\frac{\dt}{6120}\left(
\begin{array}{ccc}
 29123 & -32576 & 15789 \\
 12973 & -9456 & 6779 \\
 3963 & 1424 & 2869 \\
\end{array}
\right)
 \, \left (  \begin{array}{cccccc}
 f \left( v_{n+2/3},  t_{n+2/3} \right ) \\
f \left( v_{n+1/3},  t_{n+1/3} \right ) \\
 f \left( v_n, t_{n} \right )
\end{array} \right ),
\end{eqnarray}
and has a local truncation error of  third order,
\begin{eqnarray} \label{EIS_3_step_20b}
\vtau_{n}  &= & \frac{1}{991440}\left( \begin{array}{ccc}
115733 \\
 33623\\
 5573
\end{array} \right ) \frac{d^4} {dt^4}  u(t_n) \,\dt^3 \, + \, O(\dt^4) \nonumber \\
\nonumber \\
& & \approx \; \left( \begin{array}{ccc}
0.116732 \\
0.0339133\\
0.00562112
\end{array} \right ) \frac{d^4} {dt^4}  u(t_n) \,\dt^3 \, + \, O(\dt^4) \;.
\end{eqnarray}

The second method is
\begin{eqnarray} \label{EIS_3_step_c3}
V_{n+1}  &= &  \left(
\begin{array}{ccccc}
- \frac{101}{96} &  \frac{97}{24} & - \frac{191}{96}\\ 
- \frac{101}{96} &  \frac{97}{24} & - \frac{191}{96}\\ 
- \frac{101}{96} &  \frac{97}{24} & - \frac{191}{96}\\ 
\end{array}
\right) V_n +  \dt \left(
\begin{array}{ccc}
\frac{733}{144} &  -\frac{431}{72} & \frac{23}{12} \\ 
\frac{353}{144} & - \frac{53}{24} & \frac{4}{9}\\
\frac{47}{48} & - \frac{31}{72} & - \frac{7}{36} \\
\end{array}
\right) 
 \, \left (  \begin{array}{cccccc}
 f \left( v_{n+2/3},  t_{n+2/3} \right ) \\
f \left( v_{n+1/3},  t_{n+1/3} \right ) \\
 f \left( v_n, t_{n} \right )
\end{array} \right ). \nonumber \\
\end{eqnarray}
The  truncation error  is also third order
\begin{eqnarray} \label{EIS_3_step_20}
\vtau_{n}  &= &  \left( \begin{array}{ccc}
\frac{5303}{46656} \\
\frac{1439}{46656}\\
 \frac{119}{46656}\\
\end{array} \right ) \frac{d^4} {dt^4}  u(t_n) \,\dt^3 \, + \, O(\dt^4) \\
&= &  \left( \begin{array}{ccc}
{0.113662} \\ {0.0308428} \\ {0.00255058} \\
\end{array} \right ) \frac{d^4} {dt^4}  u(t_n) \,\dt^3 \, + \, O(\dt^4)  \nonumber
\end{eqnarray}

\begin{figure*}[t!]
    \centering
    \begin{subfigure}[t]{0.5\textwidth}
        \centering
        \includegraphics[width=0.95\textwidth]{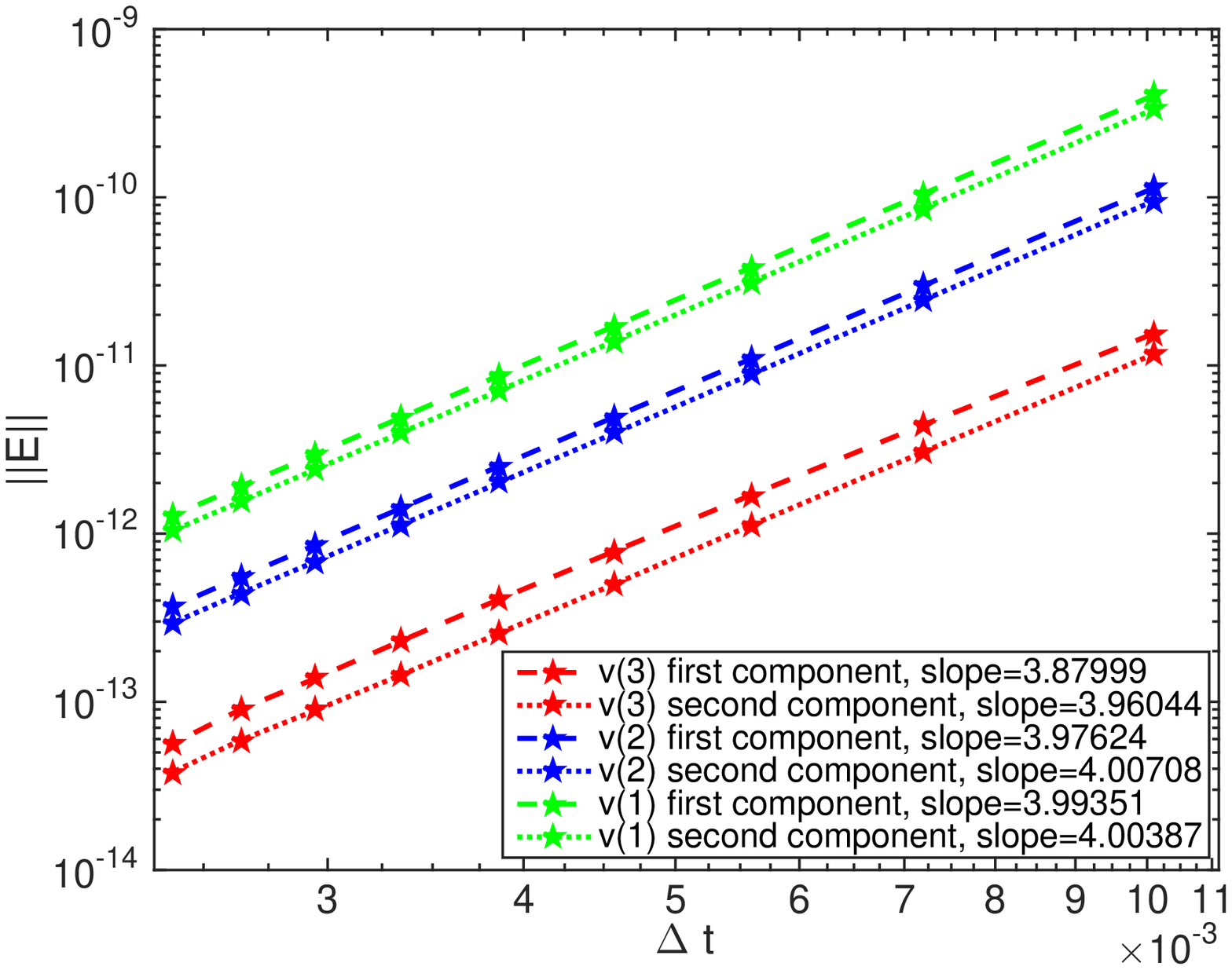}
        \caption{}
    \end{subfigure}%
    ~ 
    \begin{subfigure}[t]{0.5\textwidth}
        \centering
      \includegraphics[width=0.95\textwidth]{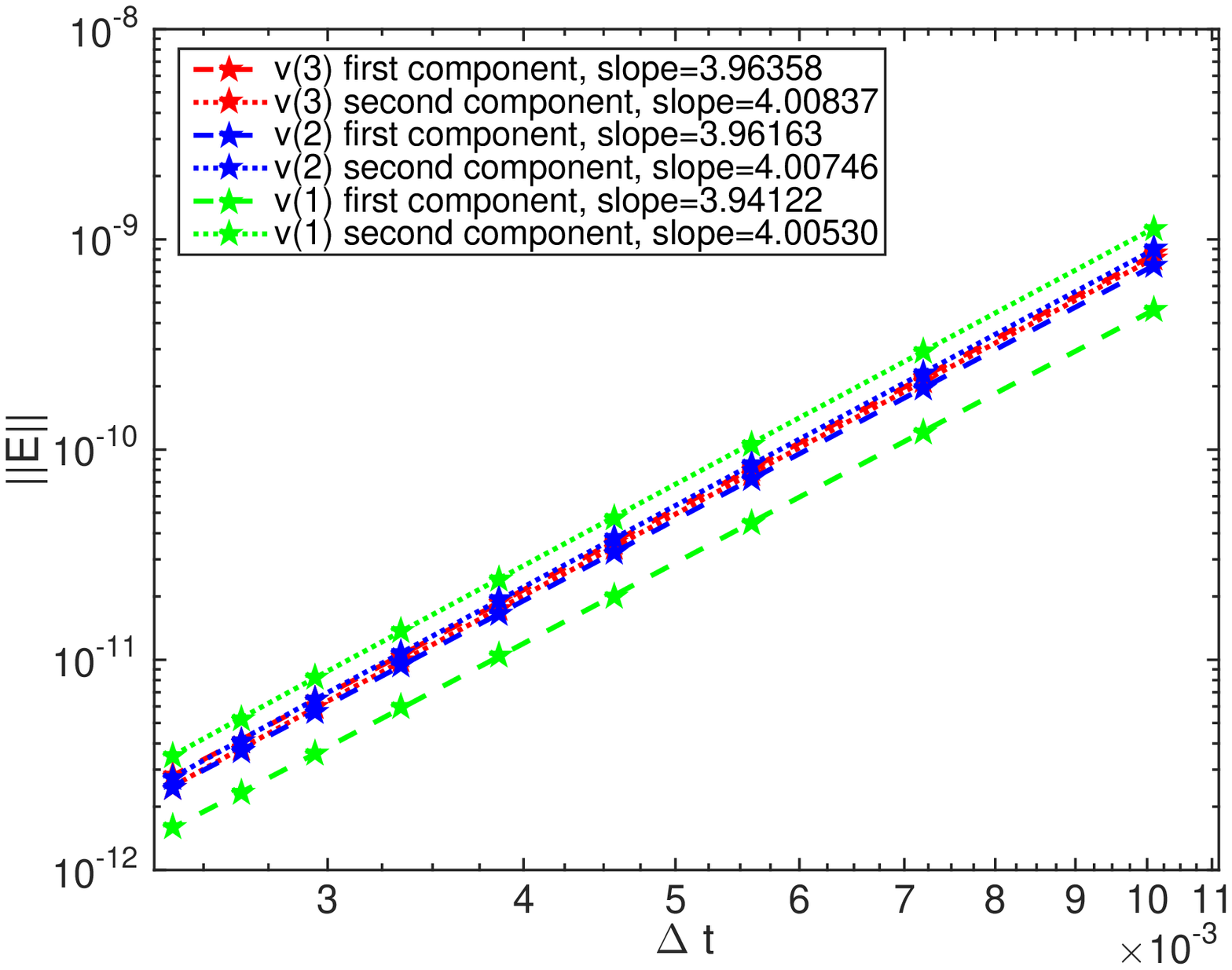}
        \caption{}
    \end{subfigure}
    \caption{Convergence plots van der Pol equation \eqref{VDP_1}. The plots show the 
    errors  vs. $ \dt$ for each component of the solution, computed for several values of $\dt$
    for  (a) the scheme \eqref{EIS_3_step_10b} and
    (b) the scheme \eqref{EIS_3_step_c3}. 
As expected, we observe fourth
order accuracy for the errors, although the truncation errors computed above are third order.}
 \label{fig:EIS_3_step_B_example}
\end{figure*}

Both these methods satisfy 
\[Q_n \vtau_{n}  = O(\dt \vtau_n) =  O(\dt^4) \;  \]
as well.  As above, this property results in an error inhibiting mechanism that produced a global error of order four. 
This can be seen once again in Figure \ref{fig:EIS_3_step_B_example}, using the nonlinear problem
 \eqref{VDP_1}  above. The results of method \eqref{EIS_3_step_10b}  are on the left and of \eqref{EIS_3_step_c3}
 are on the right.

%Both these methods have truncation errors that can be written as a linear combination of the two eigenvectors of their 
%$A$ matrix that correspond to the zero eigenvalues.
%As above, this property results in an error inhibiting mechanism that produced a global error of order four. 
%This can be seen once again in Figure \ref{fig:EIS_3_step_B_example}, using the nonlinear problem
% \eqref{VDP_1}  above. The results of method \eqref{EIS_3_step_10b}  are on the left and of \eqref{EIS_3_step_c3}
% are on the right.

\section{Conclusions} \label{conclusions} \vspace{-.1in}
While it is generally assumed that the global error will be of the order of the local truncation error,
in this work we presented an approach to creating methods that have a global error of higher order than 
predicted by the local truncation error. To accomplish this, we used the block formulation of a method
$ V_{n+1} = Q_n V_n $
where the discrete solution operator $Q_n = A+\dt B F_n $
is comprised of matrices of coefficients $A$ and $B$, and the matrix operator $F_n.$

We show that if $A$ is a diagonalizable matrix of rank one,
that has only one nonzero eigenvalue,  $z_1=1$, that corresponds to the eigenvector of all ones,
then the error inhibiting property will occur if the leading part of the local truncation error error for the linear 
constant coefficient case ($F_n=F=$ a constant)  is spanned by the  eigenvectors corresponding to the zero eigenvalues of $A$ (to the leading order).
We show that a method that has these properties will have a global error that has higher order  than the local error,
on nonlinear problems.

After presenting the concept behind these methods we use the theoretical properties above to develop
block one-step methods that are in the family of Type 3 DIMSIM methods presented in \cite{butcher1993a}.
We demonstrate in numerical examples on nonlinear problems (including a nonlinear system) that these
methods have global error that is one order higher than the local truncation errors. We also show that this is
in contrast to another Type 3 DIMSIM method which has a matrix $A$ that satisfies the first three properties {\bf C1 -- C3}, 
but does not satisfy the error inhibiting property  {\bf C4}, that the local truncation error is in the space spanned by the eigenvectors of $A$ that
correspond to the zero eigenvalues, and indeed does not give us a global error that is higher than the 
local truncation error on nonlinear test problems.

The major development in this work is the concept of an error inhibiting method and the new approach for developing
methods that are constructed to control the growth of the local truncation error. While the 
newly developed methods presented in this work can be used in place of currently standard methods (particularly in place of 
type 3 DIMSIM methods) to obtain  higher order accuracy, it is not yet known how they compare to other methods in terms of 
other important properties. In future work we intend to  the study  of the computational efficiency and storage requirements 
of these methods and the analysis of their linear stability regions. We expect that this will also lead to further development 
of error inhibiting methods that have other favorable properties.

{\bf Acknowledgements:} {\em The authors wish to thank Professor John Butcher for a very helpful discussion, and 
in particular for his valuable advice on general linear methods, especially the Type 3 DIMSIM methods. \\
The work of Sigal Gottlieb was supported by AFOSR grant FA9550-15-1-0235.}

\bibliographystyle{amsplain}

\bibliography{EIS_ODE_ref}

\end{document}